\documentclass{amsart}
\usepackage{amsmath,amsfonts,amssymb,amsthm,xcolor,amscd,latexsym,xspace,enumerate,bigints}
\usepackage[pagebackref, colorlinks=true,linkcolor=blue,citecolor=red]{hyperref}
\usepackage{pdfpages}
\usepackage{tikz-cd}
\usetikzlibrary{arrows}
\usepackage{graphicx}
\usepackage{mathtools}
\usepackage{array,booktabs}
\usepackage{booktabs,tabularx}
\usepackage{wrapfig}
\usepackage{graphicx}
\usepackage{changepage}
\usepackage{bigints}
\usepackage{comment}
\usepackage[margin=1in]{geometry}
\usepackage{tikz}

\graphicspath{ {./images/} }
\makeindex

\newcommand{\R}{\mathbb{R}}
\newcommand{\Z}{\mathbb{Z}}

\newcommand{\D}{\mathcal{D}}

\newcommand{\ar}{\longrightarrow}

\DeclareMathOperator{\Hom}{Hom}
\DeclareMathOperator{\rank}{rank}
\DeclareMathOperator{\dgm}{dgm}

\newtheorem{theorem}{Theorem}[section]

\newtheorem{lemma}[theorem]{Lemma}
\newtheorem{proposition}[theorem]{Proposition}
\newtheorem{definition}[theorem]{Definition}

\newtheorem{remark}[theorem]{Remark}
\newtheorem{example}[theorem]{Example}

\usepackage{mdframed} 

\title[A result of continuous stability for generalized persistence diagrams]{A result of continuous stability\\ 
for generalized persistence diagrams}

\author[Y. Bavuma]{Yanga Bavuma}

\author[M. Mijangos]{Martin Mijangos}

\author[F.G. Russo]{Francesco G. Russo}

\thanks{\\
Y. Bavuma - University of Cape Town - Cape Town - South Africa\\
M. Mijangos - Universidad Nacional Aut\'onoma de M\'exico - Ciudad de M\'exico - M\'exico\\
F.G. Russo - Universit\'a degli Studi di Camerino - Camerino - Italia}

\begin{document}

\begin{abstract}
We study persistence modules obtained from the $A$-fundamental groupoids of finite graphs. Moreover, we verify the categorical properties needed to form generalized persistence diagrams, and describe a quotient that identifies equivalent groupoids. A normalized linear presheaf construction places these modules in an essentially small abelian category, where Patel's Continuity Theorem gives a $1$-Lipschitz estimate in erosion distance. Explicit graph filtrations illustrate the groupoid-valued construction, and so we suggest a new approach to interpret persistence diagrams and large data sets in topological data analysis,  based on recent works of  A-homotopy theory of the last years.
\end{abstract}
\maketitle

\section{Introduction and statement of the main results}

Finite graphs support several complementary kinds of algebraic information. In the study of Cayley graphs of Pauli groups, Bavuma, D'Angeli, Donno and Russo \cite{badadoru} relate central products of groups to graph constructions and use this structure to compute adjacency spectra. Bavuma, Russo and Spessato \cite[Theorems~4, 6--8 and Section~4.2]{barusp} develop a corresponding approach to subdirect products through graph bundles and pullbacks. In particular, their network $K$-theory assigns a Grothendieck group to a graph and a fixed fibre; its dependence on the base graph is contravariant \cite[Definition~9 and Corollary~3]{barusp}. These contributions motivate the systematic study of algebraic structures attached to graphs and of their behaviour under graph operations and morphisms. A complementary motivation for the present paper comes from data analysis: Bravetti, Mijangos and Padilla-Longoria \cite{martinealessandro} use persistent homology of weighted transition graphs to compare musical styles. Their approach illustrates why a graph invariant should be studied across a filtration, together with a quantitative notion of stability. It also suggests considering other algebraic targets when the information of interest is not adequately represented by vector-space homology alone. The spectral viewpoint has a classical background in \cite{spectrabook}; the earlier network-oriented development of combinatorial homotopy appears in \cite{psapm,foundations}. For this purpose, $A$-homotopy theory provides an established homotopy theory built directly from graph maps and discrete paths \cite[Section~2]{Babson2006}. Kapulkin and Mavinkurve \cite[Definition~3.4 and Proposition~3.5]{kapulkin_grupoid} construct the $A$-fundamental groupoid and its functoriality. The groupoid retains the vertices, the connected components and the homotopy classes of paths between vertices, without requiring a preferred base vertex in every component. Their covering theory and Seifert--van Kampen Theorems \cite[Sections~4--5]{kapulkin_grupoid} further demonstrate the usefulness of this categorical viewpoint. The construction of this groupoid is therefore part of the background of the present paper.

The persistence framework is supplied by Patel \cite[Sections~2, 6--8]{patel2016generalized}. His type-$\mathcal A$ construction uses images in an essentially small symmetric monoidal category, whereas type-$\mathcal B$ diagrams use the additional relations supplied by short exact sequences in an abelian category. The latter satisfy a continuous stability estimate \cite[Theorem~8.2]{patel2016generalized}. The distinction matters here because a category of groupoids is not itself abelian. The broader interest in persistence beyond vector spaces is also reflected in the ``saecular approach'' of Ghrist and Henselman-Petrusek \cite[Introduction]{ghristhenspet2021}.

Taken together, these perspectives motivate a common question: how can one associate persistence diagrams to graph filtrations while retaining a groupoid-valued stage and controlling the effect of perturbations? Our approach combines the functor of \cite{kapulkin_grupoid} with the framework of \cite{patel2016generalized}. It is broader in its choice of categorical target than the usual vector-space construction. The connections with \cite{badadoru,barusp,martinealessandro} are methodological motivations; recovering their spectral, bundle-theoretic or application-specific results would require additional comparison theorems. In particular, the contravariant network $K$-theory of \cite{barusp} is not identified here with our covariant groupoid persistence.

We state the four results together, making explicit the conventions needed to read them. All graphs are finite, simple and undirected; we also allow the empty graph. A \emph{graph map} sends adjacent vertices to adjacent or equal vertices \cite[Definitions~1.1--1.2]{kapulkin_grupoid}. A \emph{groupoid} is a category in which every morphism is invertible. Unless otherwise specified it is small. A category is \emph{essentially small} if it is equivalent to a small category; see \cite[Section~IV.4]{maclane1998categories}. A symmetric monoidal structure consists of a bifunctor $\otimes$, a unit $I$, and coherent associativity, unit and symmetry isomorphisms \cite[Sections~VII.1 and XI.1]{maclane1998categories}. A categorical image is a factorization through a monomorphism universal among such factorizations \cite[Section~2, p.~399]{patel2016generalized}.

Write $\mathbf{Gpd}_{f_0,f_g}$ for the category of small groupoids with finitely many objects and a finite set of generating morphisms, with all functors as morphisms. Generation allows composition, inversion and identities. Finite presentation additionally requires a finite list of relations; it is not imposed on the whole target category. This distinction permits closure under the image construction below.

\begin{theorem}[First Main Theorem]\label{main1}
The category $\mathbf{Gpd}_{f_0,f_g}$ is essentially small and symmetric monoidal under
\[
\mathcal G\otimes\mathcal H=\mathcal G\sqcup\mathcal H,
\qquad I=\varnothing.
\]
It has categorical images: the image of a functor $F:\mathcal G\to\mathcal H$ is the subgroupoid of $\mathcal H$ generated by its image objects and image morphisms.
\end{theorem}

The monoidal construction is the standard one furnished by coproducts \cite[Section~VII.1, p.~163; Section~XI.1]{maclane1998categories}. The assertion to check is that the stated finiteness conditions survive it and the categorical image construction. The image need not consist only of morphisms that are individual values of $F$: identifying objects can make additional composites possible.

For a graph $\Gamma$, the established groupoid $\Pi_1^A(\Gamma)$ has vertices as objects and endpoint-preserving $A$-homotopy classes of discrete paths as morphisms. Paths are finite vertex sequences with consecutive vertices adjacent or equal; homotopies use rectangular grids with the Cartesian graph product, after allowing stationary steps. Composition is concatenation and inverses reverse paths \cite[Sections~1.3--1.5, 2; Definition~3.4]{kapulkin_grupoid}.

\begin{theorem}[Second Main Theorem]\label{main2}
The $A$-fundamental groupoid functor restricts to
\[
\Pi_1^A:\mathbf{Graph}_{\mathrm{fin}}\longrightarrow
\mathbf{Gpd}_{f_0,f_g}.
\]
Its objects are the vertices of $\Gamma$, and the classes of oriented edges form a finite generating set for $\Pi_1^A(\Gamma)$.
\end{theorem}

Functoriality is that of \cite[Definition~3.4 and the following paragraph]{kapulkin_grupoid}; the finite-generation assertion follows directly from the path description. Thus a graph filtration yields a groupoid persistence module by composition. The next result records a classical criterion for an individual groupoid, as distinct from the category of groupoids in Theorem~\ref{main1}.

\begin{theorem}[Third Main Theorem]\label{main3}
Let $\mathcal G$ be a possibly large groupoid equipped with a symmetric monoidal structure $(\otimes,I,\alpha,\lambda,\rho,\sigma)$. It is essentially small if and only if its connected components admit a set of representatives and every group $\operatorname{Aut}_{\mathcal G}(X)$ is a set. Every morphism $f:X\to Y$ has categorical image $(Y,\operatorname{id}_Y)$. Moreover, $\operatorname{Aut}_{\mathcal G}(I)$ is abelian.
\end{theorem}

The size criterion is the skeleton argument of \cite[Section~IV.4, pp.~93--94]{maclane1998categories}; the last assertion is an application of the Eckmann--Hilton argument \cite[Theorem~1.3.1, pp.~10--11]{ronbrown2}. The monoidal structure is part of the hypothesis: abelianness of one vertex group is not asserted to construct such a structure.

For the fourth result fix a field $k$ of characteristic zero, put $\mathcal C=\mathbf{Gpd}_{f_0,f_g}$, and choose a small full subcategory $\mathcal D\subset\mathcal C$ containing one representative of each isomorphism class of \emph{nonempty} objects. For $G\in\mathcal C$ we consider
\[
L_k(G)(H)=k[\operatorname{Hom}_{\mathcal C}(H,G)],\qquad H\in\mathcal D,
\]
where $k[S]$ is the vector space with basis $S$. Precomposition defines a presheaf on $\mathcal D$, and postcomposition defines $L_k$ on functors. Let $\mathcal A_k$ be the smallest full subcategory of $[\mathcal D^{\mathrm{op}},\mathbf{Vect}_k]$, closed under isomorphisms, finite direct sums, kernels and cokernels, that contains these presheaves. This is a linearized Yoneda construction \cite[Section~III.2, pp.~61--62]{maclane1998categories}; the restriction to nonempty test objects ensures $L_k(\varnothing)=0$.

A persistence module is \emph{constructible} if, for some finite set of parameters, it is constant up to isomorphism between consecutive parameters and is the monoidal unit before the first one \cite[Definitions~2.1--2.2]{patel2016generalized}. Denote by $d_I$ the infimum of shifts admitting an interleaving (natural maps in both directions whose composites are the structure maps for twice the shift), and by $d_E$ the erosion distance between diagrams \cite[Definitions~3.1--3.2 and 5.1--5.3]{patel2016generalized}.

\begin{theorem}[Fourth Main Theorem]\label{main4}
The category $\mathcal A_k$ is essentially small, concrete, $k$-linear and abelian, and $L_k:\mathcal C\to\mathcal A_k$ is faithful with $L_k(\varnothing)=0$. Moreover, $\mathcal A_k$ is not equivalent to $\mathbf{Ab}$, the category of abelian groups.
For every constructible $M:(\mathbb R,\leq)\to\mathcal C$, the module $\widehat M=L_k\circ M$ is constructible and has a type-$\mathcal B$ diagram, valued in the Grothendieck group of $\mathcal A_k$ modulo short exact sequences. For any two such modules,
\[
d_E(\widehat M_{\mathcal B},\widehat N_{\mathcal B})
\leq d_I(\widehat M,\widehat N)
\leq d_I(M,N).
\]
\end{theorem}

The first inequality is Patel's Continuity Inequality \cite[Theorem~8.2, p.~416]{patel2016generalized}, once the target and constructibility have been verified. We use the effective-class preorder on the Grothendieck group; the proof uses only its inequalities, as explained below. The second inequality follows because postcomposition preserves interleavings. Faithfulness concerns the functor $L_k$; it does not imply that the resulting Grothendieck-group diagram is a complete invariant of the original groupoids.

The paper is organized in the following way: after the statement of the main results in the present introduction, Sections \ref{s2}, \ref{s3} and \ref{sec:persistence-modulo-equivalence} present the minimal notions which we required from the references, in order to have an appropriate uniform feedback on category theory, persistence diagrams and stability in categories which are appropriate to study the topological data analysis. The proofs of the first two main results are both placed  in Section \ref{s5}, while Section \ref{s6} is devoted to the proof of our third main result. Here an argument of Eckmann and Hilton is adpated to the present context of study. Finally, Section \ref{sec:linearization} involves the proof of our last main result and adapts the so-called Yoneda Embedding Theorem to the case of those categories of graphs which turn out to be crucial in our investigations. The terminology and the notations follow mostly those in \cite{ronbrown1, ronbrown2,  maclane1998categories} for the part on category theory and homotopy theory, while the inspiration from the techniques of topological data analysis, which are at the basis of our original motivations, is taken from methodology and ideas in \cite{gc, edelsbrunner, patel2016generalized}.

\section{Preliminaries of categorical nature and \texorpdfstring{$A$}{A}-homotopy groups}\label{s2}

In this paper, we will deal with the categorical definitions of groups and groupoids, and thus will need to describe the properties we will need for our main results in this categorical language. The algebraic intuition of these concepts is enough to understand the overall paper, but the categorical language is more powerful for our framework. First, we will start with the usual definition of a group as a category.

\begin{definition}[{\cite[Sections~I.2 and I.5, pp.~11 and 20]{maclane1998categories}; \cite[Section~6.1]{ronbrown1}}]
A \textit{group} $G$ is a category with one object, and all its morphisms are isomorphisms.    A \textit{groupoid} $\mathcal{G}$ is any (small) category where all morphisms are isomorphisms.    
\end{definition}

As we can see, for a groupoid we drop the condition that it needs to have only one object. Meaning, a group is the special case of a groupoid where it only has one object. We use presentations by generating arrows and relations, as in \cite[Sections~8.2--8.3]{ronbrown1}. A finitely presented groupoid has a finite object set, finitely many generating arrows and finitely many relations between parallel words. We do not use the distinct categorical definition of a finitely presentable object via preservation of filtered colimits.

\begin{definition}[{Pair groupoid; \cite[Section~1.7, p.~26]{ronbrown2}}]
Let \(S\) be a set. The \textit{pair groupoid} on \(S\), denoted by
\(\operatorname{Pair}(S)\), is the groupoid whose set of objects is \(S\)
and such that, for every \(x,y\in S\), there is exactly one morphism
\[
u_{xy}\colon x\longrightarrow y.
\]
Composition, identities, and inverses are necessarily given by
\[
u_{yz}\circ u_{xy}=u_{xz},
\qquad
\operatorname{id}_x=u_{xx},
\qquad
u_{xy}^{-1}=u_{yx}.
\]
For \(n\geq 1\), we write
\[
\operatorname{Pair}_n
=
\operatorname{Pair}(\{1,\ldots,n\}).
\]

If \(G\) is a group, we denote by \(BG\) the one-object groupoid whose
unique object is denoted by \(*\) and whose automorphism group is \(G\);
that is,
\[
\operatorname{Ob}(BG)=\{*\},
\qquad
\Hom_{BG}(*,*)=G,
\]
with composition given by the group operation.
\end{definition}

Throughout this section, we consider a fixed field $\mathbb{F}$.
For vector-space persistence we follow \cite[Section~2, pp.~195--196]{chazal2014persistence_geom}.

\begin{definition}
A persistence module $V$ over the real numbers \cite[Section~2]{chazal2014persistence_geom} is a family of vector spaces $\{V_a|a\in \R\}$ together with a family of linear maps $\{v_a^b:V_a\ar V_b\}$ such that $v_b^c\circ v_a^b=v_a^c$ for all $a\leq b\leq c$ and $v_a^a=Id_{V_a}$. We will denote by $(V, v_a)$ a persistence module or just by $V$ if the persistence structure is clear. 
\end{definition}

One way to obtain a persistence module in data analysis is to take the simplicial homology of filtered complexes. More precisely, a \textit{filtration} of  a simplicial complex $K$ is a family of simplicial 
complexes $\{K_a\}_{a\in \R}$ such that $K_a\subseteq K_b$ whenever $a\leq b$. For $p$ a non-negative integer, we take the $p$-simplicial homology of the complexes $K_a$ with coefficients in $\mathbb{F}$, $H_p(K_a)$. For any $a\leq b $ the inclusions induce maps $v_a^b:H_p(K_a)\ar H_p(K_b)$ which for functoriality properties satisfy $v_b^c\circ v_a^b=v_a^c$ for all $a\leq b\leq c$. In this way, $\{H_p(K_a)\}_{a\in \R}$ defines a persistence module. Given a pair of persistence modules, we have a measure of dissimilarity between them defined as follows. 
\begin{definition}
    Let $(U, u_a)$ and $(V, v_a)$ be persistence modules and $\epsilon$ be a real number. A homomorphism $\phi:U\ar V$ of degree $\epsilon$ \cite[Section~2]{chazal2014persistence_geom} is a collection of maps $$\{\phi_a:U_a\ar V_{a+\epsilon}\}_{a\in \R}$$
    such that $v_{a+\epsilon}^{b+\epsilon}\circ\phi_a=\phi_{b}\circ u_a^b$ for all $a\leq b$. We denote by $\Hom ^\epsilon(U, V)$ the set of all these homomorphisms.
\end{definition}

For $\epsilon \geq 0$ we can take, in particular, the family of maps $v_a^{a+\epsilon}$, which is an element of $\Hom^\epsilon(V, V)$ and is denoted by $1_V^\epsilon$.
\begin{definition}
    For $\epsilon\geq0$, two persistence modules $U$ and $V$ are $\epsilon$-interleaved if there are maps $\phi\in \Hom^\epsilon(U,V)$ and $\psi\in\Hom^\epsilon(V,U)$ such that $\phi\circ\psi=1_V^{2\epsilon}$ and $\psi\circ\phi=1_U^{2\epsilon}$, with the shifts understood in composition \cite[Section~2]{chazal2014persistence_geom}.
\end{definition}

It turns out that, under some assumptions about the persistence modules, these can be summarized in the so-called persistence diagrams.

\begin{definition}
    A finite persistence diagram is the union of a finite multiset of points $(x, y)\in \bar{\R}^2$ where $x<y$,
    with the set of points of the form $(x, x)\in \bar{\R}^2$, taken with infinite multiplicity.
\end{definition}

We denote by $\D$ the set of these finite diagrams. Diagrams of general $q$-tame modules need not be finite; the stability theorem below uses the more general notion in \cite[Section~2]{chazal2014persistence_geom}. In order to ensure the stability of the association of a persistence diagram, we endow $\D$ with a metric. 

\begin{definition}
    Let $1\leq p<\infty$ and $A, B\in \D$. The \emph{$p$-Wasserstein distance} $d_p$ is defined as
    \begin{equation}\label{eq:wasserstein}
    d_p(A, B)=\left(\inf_\gamma\sum_{a\in A}||a-\gamma(a)||_\infty^p \right)^{1/p}
    \end{equation}
    where the infimum runs over all the bijections $\gamma:A\ar B$ and
    we define as usual $$ ||a-b||_\infty=\max\{|x_1-x_2|, |y_1-y_2|\}$$
    for $a=(x_1, y_1)$ and $b=(x_2, y_2)$,
    setting $|y_1-y_2|$ to zero if both $y_1$ and $y_2$ are infinite.
\end{definition}   

If we take the limit $p\ar \infty$ we get 
$$d_\infty(A, B)=\inf_{\gamma}\sup_{a\in A}\{||a-\gamma(a)||_\infty\}.  $$
which is usually called the \textit{bottleneck distance} and is simply denoted as $d_b$.

Recall that a persistence module is $q$-tame if $\rank u_a^b<\infty$ whenever $a<b$ \cite[Section~2, Theorem~2.3]{chazal2014persistence_geom}. The strict inequality is essential: finite-dimensional values at $a=b$ are not required. 

\begin{theorem}[{\cite[Theorem~2.3]{chazal2014persistence_geom}}]
If $(U, u_a)$ is a q-tame persistence module, then it has associated a persistence diagram denoted as $\dgm(U)$. If $U$ and $V$ are q-tame $\epsilon$-interleaved persistence modules, then the bottleneck distance between the diagrams is less or equal than $\epsilon$.
\end{theorem}

Discrete homotopy theory (or $A$-homotopy theory) allows us to study the combinatorial ``holes'' of graphs directly from their graph structure by using discrete grids of a graph instead of continuously deforming loops in topological spaces. Following \cite[Definition~1.1]{kapulkin_grupoid}, a graph $\Gamma = (V, E)$ consists of a set of vertices $V$ and a set of edges $E \subseteq \{\{u,v\} \ | \ u,v \in V, \ \ u \neq v\}$, comprising unordered pairs of distinct vertices. We write $u \sim v$ to denote $\{u, v\} \in E$. It is also useful to recall, according to \cite[Example~1.3(1)]{kapulkin_grupoid},  that for any integer $m \ge 0$, the \textit{discrete interval graph} of length $m$, denoted $I_m$, is the path graph with the vertex set $V(I_m) = \{0, 1, \dots, m\}$ and the edge set $E(I_m) = \{\{i, i+1\} \mid 0 \le i \le m-1\}$. 

In continuous topology, a (trivial) path can ``stand still'' at a point, but the domain of the path is never trivial, i.e. it is always the unit interval $[0,1]$. To emulate this in discrete homotopy theory, we must allow our maps to map adjacent vertices to the exact same vertex, so that the graph map maintains the preservation of mapping edges to edges when we use the discrete interval graph as the domain for a path.

The notion of \textit{graph map} serves the scope of the present context: according to \cite[Definition~1.2]{kapulkin_grupoid} or to similar conditions in \cite{barusp}, if $\Gamma_1=(V_1,E_1)$ and $\Gamma_2=(V_2,E_2)$ are graphs, we call  \textit{graph map} (also called  \textit{simplicial graph homomorphism} in \cite[Definition~2.1(2)]{Babson2006})  a function $f : V_1 \to V_2$ such that for all $u, v \in V_1$ such that $u \sim v$ implies $f(u) \sim f(v) \text{ or } f(u) = f(v).$

Also the notion of Cartesian  product of graphs is quite standard in the present situation: according to \cite[Definition~1.6]{kapulkin_grupoid} and \cite{badadoru, barusp} for instance, 
given two graphs $\Gamma_1$ and $\Gamma_2$, their \textit{Cartesian product} $\Gamma_1 \square \Gamma_2$ is the graph with vertex set $V_1 \times V_2$. Two vertices $(u_1, v_1) \sim (u_2, v_2)$ are adjacent in $\Gamma_1 \square \Gamma_2$ if and only if exactly one of the following two conditions hold:
\begin{enumerate}
    \item[{\rm (i).}] $u_1 = u_2$ and $v_1 \sim v_2$ in $\Gamma_2$,
    \item[{\rm (ii).}] $u_1 \sim u_2$ in $\Gamma_1$ and $v_1 = v_2$.
\end{enumerate}

Therefore one can work with discrete paths and  based loop now. 
\begin{definition}[See \cite{kapulkin_grupoid}, Definitions 1.9, 1.25, 1.26]
A \textit{discrete path} in a graph $\Gamma$ is a graph homomorphism $p : I_m \to \Gamma$ for some $m \ge 0$.  Let $v_0 \in V$ be a fixed base vertex. A \textit{based loop} of length $m$ at $v_0$ is a discrete path $p : I_m \to \Gamma$ such that $p(0) = p(m) = v_0$. 
Two based loops $p, q \colon I_m \to \Gamma$ of the same length are \textit{$A$-homotopic}, denoted $p \simeq_A q$, if there exists an integer $n \ge 0$ and a graph homomorphism $F \colon I_m \square I_n \to \Gamma$ such that:
\begin{enumerate}
    \item[{\rm (i).}] $F(i, 0) = p(i)$ for all $0 \le i \le m$,
    \item[{\rm (ii).}] $F(i, n) = q(i)$ for all $0 \le i \le m$,
    \item[{\rm (iii).}] $F(0, j) = F(m, j) = v_0$ for all $0 \le j \le n$.
\end{enumerate}
The map $F$ is called an \textit{$A$-homotopy} from $p$ to $q$.
\end{definition}

To compare loops of different lengths $m$ and $k$, we can always extend a loop to a larger length $M \ge \max(m, k)$ by defining it to remain stationary at the base vertex $v_0$ for all steps beyond its original length. Because graph homomorphisms allow $f(x) = f(x+1)$, this extension is a valid graph map. We say that two loops of arbitrary lengths are $A$-homotopic if their extensions to some common length $M$ are $A$-homotopic. The equivalent finite-grid and shrinking-map formulations are compared in \cite[Section~2, especially Theorem~2.14]{kapulkin_grupoid}.

\begin{definition}[{The $A$-homotopy group; \cite[Definition~2.4]{Babson2006}}]
Let $(\Gamma, v_0)$ be a pointed graph, i.e. a graph with a chosen base vertex $v_0$. The \textit{$A$-homotopy group}, denoted $A_1(\Gamma, v_0)$ (the subscript $1$ is because it is the \textit{first} $A$-homotopy group), is defined as the set of $A$-homotopy equivalence classes of based loops at $v_0$:
\[
A_1(\Gamma, v_0) = \{ [p] \mid p \text{ is a based loop at } v_0 \}
\]
\end{definition}

The group operation is defined by \textit{concatenation}. Given two classes $[p]$ and $[q]$ represented by loops $p \colon I_m \to \Gamma$ and $q \colon I_k \to \Gamma$, their product $[p] * [q]$ is the equivalence class of the concatenated loop $p * q \colon I_{m+k} \to \Gamma$ defined by:
\[
(p * q)(i) = 
\begin{cases} 
p(i) & \text{if } 0 \le i \le m \\
q(i - m) & \text{if } m \le i \le m+k
\end{cases}
\]
This operation is well-defined and associative; see \cite[Definition~2.4 and the following discussion]{Babson2006} and \cite[Definition~1.29]{kapulkin_grupoid}. The identity element is the equivalence class of the constant loop of length $0$ at $v_0$. The inverse of a loop $p \colon I_m \to \Gamma$ is the reverse loop $p^{-1} \colon I_m \to \Gamma$ defined by $p^{-1}(i) = p(m - i)$. Under this operation, $A_1(\Gamma, v_0)$ forms a valid group.

\medskip

The group is functorial on \emph{pointed} graphs and basepoint-preserving graph maps:
\[
A_1:\mathbf{Graph}_{\mathrm{fin},*}\longrightarrow\mathbf{Grp};
\]
see \cite[Section~1.5]{kapulkin_grupoid}. For a connected graph, changing the base vertex gives an isomorphic group, with the isomorphism depending on a connecting path; compare \cite[Propositions~6.3.1--6.3.2]{ronbrown1}. Accordingly, $A_1(\Gamma)$ without a basepoint denotes only its isomorphism type.

\begin{proposition}
    For a connected graph $\Gamma$, we have that $\pi_1(\Gamma)$ is a free group on  finitely many generators.
\end{proposition}

\begin{proof}
    This follows from \cite[Proposition~1A.2, p.~84]{hatcher}, by choosing a maximal tree.
\end{proof}

\medskip

The topological fundamental group of a $3$-cycle or $4$-cycle, regarded as a one-dimensional cell complex, is $\mathbb Z$, whereas its $A$-homotopy group is trivial \cite[Example~1.23]{kapulkin_grupoid}. Attaching $2$-cells along these short cycles gives the following comparison with the topological fundamental group. For the required background on cell complexes, see \cite{hatcher}. This is why we omit proofs of the followng two lemmas.

\begin{lemma}[See \cite{hatcher}, Proposition 1.26(a)] \label{2cellspace}
Let $X$ be a path-connected topological space and $x_0 \in X$. If $Y$ is obtained from $X$ by attaching $2$-cells along loops of $X$, then the inclusion $X \hookrightarrow Y$ induces a surjection $$\pi_1(X, x_0) \to \pi_1(Y, x_0)$$ whose kernel is $N$, the normal subgroup generated by the attaching loops transported to $x_0$ along chosen paths. Thus $$\pi_1(Y) \approx \pi_1(X)/N.$$
\end{lemma}

The following result can be also found in  \cite[p.~32]{Babson2006}.

\begin{lemma}[See \cite{foundations}, Proposition 5.12] \label{2cellgraph}
Let $X_\Gamma$ be the topological space obtained from the graph $\Gamma$ by attaching $2$-cells to $\Gamma$, viewed as a one dimensional cell complex, one along each $3$-cycle and each $4$-cycle of $\Gamma$. Then $$A_1(\Gamma,v_0) \cong \pi_1(X_\Gamma,v_0).$$
\end{lemma}

It can be useful to argue from the perspective of A-homotopy theory in the next result, so we may start to note the differences with the usual arguments that one has in classical algebraic topology.

\begin{proposition} \label{finpres}
Let $\Gamma = (V, E)$ be a graph, and $v_0 \in V$ be a base vertex. Then its $A$-homotopy group, denoted $A_1(\Gamma, v_0)$, is a finitely presented group.
\end{proposition}

\begin{proof}
By Lemma \ref{2cellgraph}, the $A$-homotopy group $A_1(\Gamma, v_0)$ is isomorphic to the topological fundamental group $\pi_1(X_\Gamma)$ of a CW complex $X_\Gamma$ constructed by attaching $2$-cells to the $3$ and $4$-cycles of $\Gamma$. Because $\Gamma$ is a finite graph, it contains a finite number of edges and cycles, yielding a finite CW complex. By Lemma \ref{2cellspace}, the fundamental group of a finite CW complex is finitely presented, i.e. finite generators and finite relations from the finite number of loops in $\Gamma$; therefore, $A_1(\Gamma, v_0)$ is finitely presented.
\end{proof}

\section{\texorpdfstring{$A$}{A}-homotopy groupoids of graphs and persistence diagrams}\label{s3}

The following definition and its functoriality are due to the established theory in \cite[Definition~3.4]{kapulkin_grupoid}.

\begin{definition}[$A$-fundamental groupoid]\label{def:agroupoid}
For a graph $\Gamma$, the groupoid $\Pi_1^A(\Gamma)$ has object set $V(\Gamma)$. A morphism $u\to v$ is an endpoint-preserving $A$-homotopy class of finite discrete paths from $u$ to $v$, allowing stationary extensions when lengths differ. Homotopies use $I_m\square I_n$, with the endpoints fixed throughout. If $p:u\to v$ and $q:v\to w$, then $[q]\circ[p]=[p*q]$. Identities are constant paths and $[p]^{-1}$ is represented by the reversed path. These operations are well-defined and satisfy the groupoid axioms on homotopy classes.
\end{definition}

A graph map $f:\Gamma\to\Lambda$ induces the functor $\Pi_1^A(f)$ given by $u\mapsto f(u)$ and $[p]\mapsto[f\circ p]$. Moreover,
\[
\operatorname{Aut}_{\Pi_1^A(\Gamma)}(v_0)\cong A_1(\Gamma,v_0)
\]
by \cite[Proposition~3.5]{kapulkin_grupoid}. This permits all components to be treated at once. The Cartesian graph product $\square$ used for homotopies is the box product of \cite[Definition~1.6]{kapulkin_grupoid}; it is not the categorical product of the graph category, nor the disjoint-union tensor of Theorem~\ref{main1}.

\begin{proposition}
  Let $\Gamma$ be a graph. Then its $A$-homotopy groupoid is finitely presented.  
\end{proposition}

\begin{proof}
    By Proposition~\ref{finpres}, the vertex group of each component is finitely presented. A connected groupoid with object set $S$ and chosen vertex group $G$ is noncanonically isomorphic to $\operatorname{Pair}(S)\times BG$; see \cite[Section~1.7, equation~(1.7.1), p.~26]{ronbrown2}. For finite $S$, finitely many connecting arrows, together with generators and relations for $G$, give a finite groupoid presentation. There are only finitely many components.
\end{proof}

To utilize the generalized persistence framework developed by  Patel \cite{patel2016generalized}, we must define our persistence modules over categories equipped with specific algebraic structures \cite{maclane1998categories}. A \textit{symmetric monoidal category}  $(\mathcal{C}, \otimes, e)$ (according to \cite[Sections~VII.1 and XI.1]{maclane1998categories}) is a category $\mathcal{C}$ equipped with a bifunctor $\otimes \colon \mathcal{C} \times \mathcal{C} \to \mathcal{C}$ (the tensor product), an identity object $e \in \mathcal{C}$ and natural isomorphisms that satisfy constraints of associativity, unitality, and symmetry. As it often happens in $K$-theory, one can introduce an appropriate notion of Grothendieck group. 

\begin{definition}
Let $(\mathcal{C}, \otimes, e)$ be an essentially small symmetric monoidal category. The Grothendieck group of $\mathcal{C}$ (in the sense of \cite[Definition~6.1.1, p.~406]{patel2016generalized}) is denoted  by $K_0(\mathcal{C})$ and is the abelian group generated by the isomorphism classes $[X]$ of objects $X \in \mathcal{C}$, subject to the relations $[X]+[Y]=[X\otimes Y]$. This is Patel's $\mathcal A(\mathcal C)$, which we denote by $K_0(\mathcal C)$; it is distinct from the quotient by short exact sequences used below.
\end{definition}

Now we may introduce the concepts which are typical of topological data analysis and inspired the origin of the present investigations.

\begin{definition}
Let $\mathcal{C}$ be an essentially small symmetric monoidal category with images. A persistence module (in the sense of Patel \cite[Definition~2.1]{patel2016generalized}) is a functor $F \colon (\mathbb{R}, \leq) \to \mathcal{C}$ out of the poset of real numbers.
Furthermore, let $S = \{s_1 < s_2 < \dots < s_n\}$ be a finite set of real numbers, and let $e \in \mathcal{C}$ be the monoidal identity object. The persistence module $F$ is said to be \textit{$S$-constructible} (in the sense of Patel \cite[Definition~2.2]{patel2016generalized}) if:
\begin{itemize}
    \item for $p \leq q < s_1$, $F(p \leq q)$ is the identity on $e$;
    \item for $s_i \leq p \leq q < s_{i+1}$, $F(p \leq q)$ is an isomorphism;
    \item for $s_n \leq p \leq q$, $F(p \leq q)$ is an isomorphism.
\end{itemize}
A persistence module $F$ is called \textit{constructible} if there exists a finite set $S$ such that $F$ is $S$-constructible.
\end{definition}

Let $\mathrm{PMod}(\mathcal{C})$ denote the full subcategory of the functor category $[(\mathbb{R}, \le), \mathcal{C}]$ consisting of constructible persistence modules. To define the persistence diagram of such modules, we follow Patel's formalization over the poset of intervals, setting the target abelian group to $\mathcal{G} = K_0(\mathcal{C})$.
Following \cite[Definition~4.1]{patel2016generalized}, define $(\mathsf{Dgm}, \supseteq)$ as the poset of all half-open intervals $[q, r) \subset \mathbb{R}$, for $q < r$, and all half-infinite intervals $[q, \infty) \subset \mathbb{R}$. The poset relation is the containment relation.

\begin{definition}
Let $S = \{s_1 < \dots < s_n\}$ be a finite nonempty set of real numbers and $\mathcal{G}$ an abelian group. We use the following endpoint convention for \textit{$S$-constructible rank maps} (in the sense of Patel \cite[Section~7]{patel2016generalized}), adapted to right-continuous modules and the left-limit ranks. Put $s_{n+1}=\infty$. Require $X([p,q))=0$ if $p<s_1$. Otherwise let $i=\max\{h:s_h\leq p\}$ and let $j$ be the least index with $q\leq s_j$, allowing $j=n+1$. Require
\[
X([p,q))=X([s_i,s_j)).
\]
Thus the left endpoint is rounded down and the right endpoint up. This specifies the endpoint and zero-boundary conventions needed for the rank-map use of \cite[Definition~4.2 and Theorem~4.1]{patel2016generalized}. 
\end{definition}

In the present situation a map $Y \colon \mathsf{Dgm} \to \mathcal{G}$ is  called \textit{$S$-finite} (following \cite[Definition~4.3]{patel2016generalized}), if $Y(I) \neq e$ implies $I = [s_i, s_j)$ or $I = [s_i, \infty)$. We say that a map $Y \colon \mathsf{Dgm} \to \mathcal{G}$ is \textit{finite} if it is $T$-finite for some set $T$.

\begin{definition}[See \cite{patel2016generalized}, Definition~4.4] 
A \textit{persistence diagram} is a finite map $Y \colon \mathsf{Dgm} \to \mathcal{G}$.
\end{definition}

Therefore the following criterion can be formulated, in order to show the powerful notion which we have just introduced.

\begin{proposition}[M\"obius inversion]\label{mif}
For an $S$-constructible rank map $$X \colon \mathsf{Dgm} \to \mathcal{G}$$ in the preceding endpoint convention, there is a unique $S$-finite map $Y \colon \mathsf{Dgm} \to \mathcal{G}$ satisfying the Möbius inversion formula
$$ X(I) = \sum_{\substack{J \in \mathsf{Dgm}\\ J \supseteq I}} Y(J), $$
for each $I \in \mathsf{Dgm}$.
\end{proposition}
\begin{proof}
Write $X_{ij}=X([s_i,s_j))$ for $1\leq i<j\leq n+1$, and set $X_{0j}=0$. For $j\leq n$, put
\[
Y([s_i,s_j))=X_{ij}-X_{i,j+1}-X_{i-1,j}+X_{i-1,j+1},
\]
and put $Y([s_i,\infty))=X_{i,n+1}-X_{i-1,n+1}$; set $Y=0$ elsewhere. The upper-interval sum telescopes to $X_{ij}$ on the grid. The rounding convention gives the same sum for every other interval. Conversely, taking these differences of the asserted upper-interval sums recovers every value of $Y$, proving uniqueness. These are the finite-difference formulas of \cite[Theorem~4.1, pp.~403--404]{patel2016generalized}, with the boundary conditions explicit.
\end{proof}

Let $(\mathcal{G}, \preceq)$ be an abelian group with a translation-invariant preorder. That is, $\preceq$ is reflexive and transitive, and $a\preceq b$ implies $a+c\preceq b+c$. We allow preorders because the image of the object monoid in a Grothendieck group need not define an antisymmetric relation. When its positive cone is pointed, this is a partial order. Let $e\in\mathcal G$ be the additive identity.

\begin{definition}
Let $Y_1, Y_2 \colon \mathsf{Dgm} \to (\mathcal{G}, \preceq)$ be two persistence diagrams. A \textit{morphism} $Y_1 \to Y_2$ of persistence diagrams (in the sense of Patel \cite[Definition~5.1, p.~405]{patel2016generalized}) is the relation
$$ \sum_{\substack{J \in \mathsf{Dgm}\\ J \supseteq I}} Y_1(J) \preceq \sum_{\substack{J \in \mathsf{Dgm}\\ J \supseteq I}} Y_2(J), $$
for each $I \in \mathsf{Dgm}$. Let $\mathsf{PDgm}(\mathcal{G})$ be the resulting preordered collection of persistence diagrams. This is the relation of \cite[Definition~5.1]{patel2016generalized}; when the coefficient preorder is a partial order, so is the diagram relation.
\end{definition}

Fix an essentially small symmetric monoidal category $\mathcal{C}$ with images. For the zero-boundary convention, assume also that every morphism from the monoidal unit $e$ is monic; its image is then isomorphic to $e$. This holds for the empty groupoid and for the zero object in an abelian category, the two cases used here. We assign to each persistence module $F \in \mathrm{PMod}(\mathcal{C})$ a persistence diagram $F_{K_0} \in \mathsf{PDgm}(K_0(\mathcal{C}))$.

We start by constructing a map
$$ dF_{\mathcal{I}} \colon \mathsf{Dgm} \to \mathcal{I}(\mathcal{C}). $$

Recall that $\mathcal I(\mathcal C)$ is the monoid of isomorphism classes. We use the left-limit convention of \cite[Section~7, p.~409]{patel2016generalized}. For $I=[p,q)$ with $q<\infty$, choose $r$ with $p\leq r<q$ sufficiently close to $q$ that $F$ is constant up to isomorphism on $[r,q)$. Define
\[
dF_{\mathcal I}([p,q))=[\operatorname{im}F(p\leq r)].
\]
For $I=[p,\infty)$ choose $r\geq p$ beyond the last critical value and use the same formula. Constructibility makes the isomorphism class independent of these choices. This formulation ensures that the structure map exists even when $p$ is close to $q$.
Composing with the canonical group-completion map $\mathcal I(\mathcal C)\to K_0(\mathcal C)$ gives $dF_{K_0}$. The canonical map need not be injective for a general monoid.

\begin{definition}
The \textit{type $K_0$ persistence diagram} of $F$ (i.e: Patel's type-$\mathcal A$ diagram in \cite[Definition~7.1]{patel2016generalized}) is the Möbius inversion
$$ F_{K_0} \colon \mathsf{Dgm} \to K_0(\mathcal{C}) $$
of $dF_{K_0} \colon \mathsf{Dgm} \to K_0(\mathcal{C})$.
\end{definition}

For an abelian target $\mathcal C$, write
\[
\mathcal B(\mathcal C)=K_0(\mathcal C)/
\langle [B]-[A]-[C]:0\to A\to B\to C\to0\text{ is exact}\rangle.
\]
This is \cite[Definition~6.2.1]{patel2016generalized}. The type-$\mathcal B$ diagram $F_{\mathcal B}$ is obtained by applying this quotient to the rank map before M\"obius inversion \cite[Definition~7.2]{patel2016generalized}.

For completeness, consider \[d_I(F,G)=\inf\{\epsilon\geq0 \ \mbox{for which there are natural transformations} \  F\to G(\cdot+\epsilon)  \ \mbox{and} \  G\to F(\cdot+\epsilon) \]\[ \mbox{whose shifted composites are the}  \ 2\epsilon \ \mbox{structure maps}\}\] according to \cite[Definitions~3.1--3.2]{patel2016generalized}. For diagrams, let \[\operatorname{Grow}_\epsilon[p,q)=[p-\epsilon,q+\epsilon), \  \mbox{with} \ \infty+\epsilon=\infty, \quad \mbox{and} \quad \nabla_\epsilon Y=Y\circ\operatorname{Grow}_\epsilon.\] Then $d_E(Y,Z)$ is the infimum of the $\epsilon$ for which $\nabla_\epsilon Y\to Z$ and $\nabla_\epsilon Z\to Y$ in the preceding diagram order \cite[Definitions~5.2--5.3]{patel2016generalized}.

\begin{theorem}[{Patel's Continuity Theorem, \cite[Theorem~8.2, p.~416]{patel2016generalized}}]\label{stability}
Let $\mathcal C$ be an essentially small abelian category, realized concretely with concrete images. For constructible persistence modules $F,G$ in $\mathcal C$,
\[
d_E(F_{\mathcal B},G_{\mathcal B})\leq d_I(F,G).
\]
\end{theorem}
The concrete realization required for the image argument is discussed explicitly in \cite[Section~8, p.~413]{patel2016generalized}. We verify it directly for $\mathcal A_k$ in Section~\ref{sec:linearization}.
We use the effective-class preorder on $\mathcal B(\mathcal C)$: $a\preceq b$ when $b-a$ is the class of an object. The continuity proof uses only these inequalities and therefore also gives the displayed estimate when this relation is a preorder. One may pass to its antisymmetric quotient if an ordered coefficient group is desired.

\section{Generalized persistence diagrams modulo equivalence}
\label{sec:persistence-modulo-equivalence}

Patel's type-$\mathcal{A}$ persistence diagram is constructed from the
commutative monoid $\mathcal{I}(\mathcal{C})$ of isomorphism classes
of objects of a symmetric monoidal category $\mathcal{C}$. In some
situations, however, the objects of $\mathcal{C}$ carry a natural
notion of equivalence that is weaker than isomorphism. In this section,
we describe the corresponding quotient of Patel's construction \cite[Sections~6.1 and 7]{patel2016generalized}. The quotient-monoid and group-completion steps are formal algebraic constructions; the relevant choice here is which equivalence of objects to impose.

Throughout this section, let
\[
(\mathcal{C},\otimes,e)
\]
be an essentially small symmetric monoidal category with images. We
write
\[
\mathcal{A}(\mathcal{C})=K_0(\mathcal{C})
\]
for the group completion of $\mathcal{I}(\mathcal{C})$. An equivalence relation $\simeq$ on the objects of $\mathcal{C}$ is
called a \textit{monoidal equivalence relation} if it satisfies the
following conditions:
\begin{enumerate}
    \item[{\rm(i).}] if $X\cong Y$ in $\mathcal{C}$, then $X\simeq Y$;
    \item[{\rm(ii).}] if $X\simeq X'$ and $Y\simeq Y'$, then $X\otimes Y\simeq X'\otimes Y'.$
\end{enumerate}
The first condition implies that $\simeq$ induces a well-defined
equivalence relation on $\mathcal{I}(\mathcal{C})$: namely,
$[X]\sim[Y]$ if and only if  $X\simeq Y.$ We define
\[
\mathcal{I}_{\simeq}(\mathcal{C})
:=
\mathcal{I}(\mathcal{C})/{\sim},
\]
and denote the equivalence class of an object $X$ by
$\langle X\rangle$. Since $\mathcal{C}$ is essentially small,
$\mathcal{I}(\mathcal{C})$ is a set, and hence
$\mathcal{I}_{\simeq}(\mathcal{C})$ is also a set.

\begin{proposition}
The operation
\[
\langle X\rangle+\langle Y\rangle
:=
\langle X\otimes Y\rangle
\]
makes $\mathcal{I}_{\simeq}(\mathcal{C})$ into a commutative monoid
with identity $\langle e\rangle$.
\end{proposition}

\begin{proof}
The operation is well defined by the compatibility of $\simeq$ with
the monoidal product. Associativity, commutativity, and the identity
property follow from the associativity, symmetry, and unit
isomorphisms of the symmetric monoidal structure, together with the
implication $X\cong Y\Longrightarrow X\simeq Y.$
\end{proof}

The following notion is a further generalization of those we have seen until now.

\begin{definition}
The \textit{Grothendieck group modulo equivalence} of
$(\mathcal{C},\otimes,e,\simeq)$ is the group completion
\[
\mathcal{A}_{\simeq}(\mathcal{C})
:=
\operatorname{Gr}
\bigl(\mathcal{I}_{\simeq}(\mathcal{C})\bigr).
\]
Equivalently, $\mathcal{A}_{\simeq}(\mathcal{C})$ is the abelian
group generated by the symbols $\langle X\rangle$, one for each
$\simeq$-equivalence class of objects of $\mathcal{C}$, subject to
the relations
\[
\langle X\otimes Y\rangle
=
\langle X\rangle+\langle Y\rangle.
\]
\end{definition}

There is a canonical surjective monoid homomorphism
\[
\pi_{\mathcal I}\colon
\mathcal{I}(\mathcal{C})
\longrightarrow
\mathcal{I}_{\simeq}(\mathcal{C}),
\qquad
[X]\longmapsto\langle X\rangle.
\]
By the universal property of group completion, it induces a
surjective group homomorphism
\[
\pi_{\mathcal A}\colon
\mathcal{A}(\mathcal{C})
\longrightarrow
\mathcal{A}_{\simeq}(\mathcal{C}),
\qquad
[X]-[Y]\longmapsto
\langle X\rangle-\langle Y\rangle.
\]

More explicitly, if
\[
N_{\simeq}
=
\left\langle
[X]-[Y]\;\middle|\;X\simeq Y
\right\rangle
\]
denotes the subgroup of $\mathcal{A}(\mathcal{C})$ generated by all
differences of equivalent objects, then
\[
\mathcal{A}_{\simeq}(\mathcal{C})
\cong
\mathcal{A}(\mathcal{C})/N_{\simeq}.
\]
Let $F\colon(\mathbb{R},\leq)\longrightarrow\mathcal{C}$
be an $S$-constructible persistence module. Recall Patel's
image-based map
\[
dF_{\mathcal I}\colon
\mathsf{Dgm}\longrightarrow\mathcal{I}(\mathcal{C}).
\]

\begin{definition}
The \textit{rank map modulo equivalence} of $F$ is the map
\[
dF_{\mathcal{I}_{\simeq}}
:=
\pi_{\mathcal I}\circ dF_{\mathcal I}
\colon
\mathsf{Dgm}
\longrightarrow
\mathcal{I}_{\simeq}(\mathcal{C}).
\]
Thus $dF_{\mathcal I_{\simeq}}(I)$ is the equivalence class of the same image object used in $dF_{\mathcal I}(I)$, with the left-limit convention specified above.
\end{definition}

Let $\eta_{\simeq}\colon
\mathcal{I}_{\simeq}(\mathcal{C})
\longrightarrow
\mathcal{A}_{\simeq}(\mathcal{C})$
be the canonical group-completion homomorphism. Define
\[
dF_{\mathcal{A}_{\simeq}}
:=
\eta_{\simeq}\circ
dF_{\mathcal{I}_{\simeq}}
\colon
\mathsf{Dgm}
\longrightarrow
\mathcal{A}_{\simeq}(\mathcal{C}).
\]
Equivalently,
\[
dF_{\mathcal{A}_{\simeq}}
=
\pi_{\mathcal A}\circ dF_{\mathcal A}, \ \ 
\mbox{where} \ \ 
dF_{\mathcal A}\colon
\mathsf{Dgm}\longrightarrow\mathcal{A}(\mathcal{C})
\]
is Patel's original type-$\mathcal A$ rank map, see \cite{patel2016generalized}.

\begin{proposition}
The map $dF_{\mathcal{A}_{\simeq}}$ is well defined and
$S$-constructible in the endpoint convention above. Moreover, it is independent of the evaluation points and the constructibility set used in its definition.
\end{proposition}

\begin{proof}
Categorical images are unique up to isomorphism. Therefore, they are
also unique up to the equivalence relation $\simeq$. Patel's map $dF_{\mathcal I}$ \cite[Section~7]{patel2016generalized} is independent of the choices of evaluation points and the constructibility set, because the image
objects obtained from different choices are isomorphic. Applying
$\pi_{\mathcal I}$ shows that their $\simeq$-equivalence classes
are equal. Consequently, $dF_{\mathcal{I}_{\simeq}}$ and
$dF_{\mathcal{A}_{\simeq}}$ are well defined. Finally, $dF_{\mathcal A}$ is $S$-constructible, and composing an
$S$-constructible map with a group homomorphism preserves
$S$-constructibility. Hence
$dF_{\mathcal{A}_{\simeq}}$ is $S$-constructible.
\end{proof}

\begin{definition}
The \textit{type-$\mathcal{A}_{\simeq}$ persistence diagram} of $F$
is the Möbius inversion
\[
F_{\mathcal{A}_{\simeq}}
\colon
\mathsf{Dgm}
\longrightarrow
\mathcal{A}_{\simeq}(\mathcal{C})
\]
of the $S$-constructible map \, 
$dF_{\mathcal{A}_{\simeq}}
\colon
\mathsf{Dgm}
\longrightarrow
\mathcal{A}_{\simeq}(\mathcal{C}).$
\end{definition}

By the M\"obius inversion formula \cite[Theorem~4.1]{patel2016generalized},
$F_{\mathcal{A}_{\simeq}}$ is $S$-finite and satisfies
for every $I\in\mathsf{Dgm}$ \[
dF_{\mathcal{A}_{\simeq}}(I)
=
\sum_{\substack{J\in\mathsf{Dgm}\\J\supseteq I}}
F_{\mathcal{A}_{\simeq}}(J).
\]
This allows us to split naturally  $F_{\mathcal{A}_{\simeq}}$ via the following result.
\begin{proposition}
Let
 $F_{\mathcal A}\colon
\mathsf{Dgm}\longrightarrow\mathcal{A}(\mathcal{C})$
be Patel's type-$\mathcal A$ persistence diagram. Then
\[F_{\mathcal{A}_{\simeq}}
=
\pi_{\mathcal A}\circ F_{\mathcal A}.\]
In particular, the type-$\mathcal{A}_{\simeq}$ persistence diagram
is a quotient of the type-$\mathcal A$ persistence diagram.
\end{proposition}

\begin{proof}
For every $I\in\mathsf{Dgm}$, Proposition \ref{mif}  gives in the present situation
\[
dF_{\mathcal A}(I)
=
\sum_{\substack{J\in\mathsf{Dgm}\\J\supseteq I}}
F_{\mathcal A}(J).
\]
Applying the group homomorphism $\pi_{\mathcal A}$ yields
\[
dF_{\mathcal{A}_{\simeq}}(I)
=
\sum_{\substack{J\in\mathsf{Dgm}\\J\supseteq I}}
\pi_{\mathcal A}\bigl(F_{\mathcal A}(J)\bigr).
\]
Therefore,
$\pi_{\mathcal A}\circ F_{\mathcal A}$ is the Möbius inversion of
$dF_{\mathcal{A}_{\simeq}}$, and hence
$F_{\mathcal{A}_{\simeq}} = \pi_{\mathcal A}\circ F_{\mathcal A}.$
\end{proof}

The image of
$\mathcal{I}_{\simeq}(\mathcal{C})$ inside
$\mathcal{A}_{\simeq}(\mathcal{C})$ defines a natural positive cone
\[
P_{\simeq}
:=
\eta_{\simeq}
\bigl(\mathcal{I}_{\simeq}(\mathcal{C})\bigr).
\]
It induces a translation-invariant preorder by
$a\preceq_{\simeq}b
 \Longleftrightarrow 
b-a\in P_{\simeq}.$
This preorder is a partial order whenever $P_{\simeq}\cap(-P_{\simeq})=\{0\}.$ Under this assumption, the type-$\mathcal{A}_{\simeq}$ diagram may
be regarded as an object of
\[
\mathsf{PDgm}
\bigl(\mathcal{A}_{\simeq}(\mathcal{C})\bigr).
\]

\begin{proposition}
The quotient homomorphism
\[
\pi_{\mathcal A}\colon
\mathcal{A}(\mathcal{C})
\longrightarrow
\mathcal{A}_{\simeq}(\mathcal{C})
\]
is order preserving with respect to the preorders induced by the corresponding positive cones (and hence the orders when the cones are pointed).
\end{proposition}

\begin{proof}
Let $\eta:\mathcal I(\mathcal C)\to\mathcal A(\mathcal C)$ be the group-completion map. If $a\preceq b$, then $b-a\in\eta\bigl(\mathcal{I}(\mathcal{C})\bigr).$ Consequently,
$
\pi_{\mathcal A}(b)-\pi_{\mathcal A}(a)
=
\pi_{\mathcal A}(b-a)
\in
\eta_{\simeq}\bigl(\mathcal{I}_{\simeq}(\mathcal{C})\bigr),$
and therefore
$\pi_{\mathcal A}(a)
\preceq_{\simeq}
\pi_{\mathcal A}(b).$
\end{proof}

\begin{remark} \rm
Since $\pi_{\mathcal A}$ is order preserving, every morphism between
Patel type-$\mathcal A$ persistence diagrams induces a morphism
between their type-$\mathcal{A}_{\simeq}$ quotients. Any applicable type-$\mathcal A$ estimate can therefore be passed to the quotient. Applying \cite[Theorem~8.1]{patel2016generalized} still requires checking the concrete-image hypothesis stated at the start of that paper's Section~8; it does not follow from the existence of categorical images alone. The quotient construction by itself does not provide type-$\mathcal B$ continuity, which uses an abelian target.
\end{remark}

\section{Proofs of  First and  Second Main Theorems}\label{s5}

We may proceed to the proof of our first main result.

\begin{proof}[Proof of Theorem \ref{main1}]
Consider $\mathcal{C}=\mathbf{Gpd}_{f_0,f_g}.$
We prove separately that $\mathbf{Gpd}_{f_0,f_g}$ is essentially small, symmetric monoidal, and has categorical images.

First, $\mathbf{Gpd}_{f_0,f_g}$ is locally small. Indeed, for fixed groupoids $\mathcal{G},\mathcal{H}\in\mathbf{Gpd}_{f_0,f_g}$, a functor $F\colon\mathcal{G}\to\mathcal{H}$ is determined by functions
\[
F_0\colon\operatorname{Ob}(\mathcal{G})
\longrightarrow\operatorname{Ob}(\mathcal{H})
\ \ 
\mbox{and} \ \ F_1\colon\operatorname{Mor}(\mathcal{G})
\longrightarrow\operatorname{Mor}(\mathcal{H}),
\]
satisfying the usual functoriality conditions. Hence
\[
\operatorname{Hom}_{\mathbf{Gpd}_{f_0,f_g}}(\mathcal{G},\mathcal{H})
\subseteq
\operatorname{Ob}(\mathcal{H})^{\operatorname{Ob}(\mathcal{G})}
\times
\operatorname{Mor}(\mathcal{H})^{\operatorname{Mor}(\mathcal{G})},
\]
and therefore every hom-class is a set.

We now prove essential smallness, using the standard representative-category criterion \cite[Section~IV.4]{maclane1998categories}. Let $\mathcal{G}\in\mathbf{Gpd}_{f_0,f_g}$ and let $S$ be a finite set of morphisms that generate $\mathcal{G}$. Since $\operatorname{Ob}(\mathcal{G})$ is finite, every morphism of $\mathcal{G}$ can be represented by a finite composable word in the alphabet
\[
S\cup S^{-1}\cup
\{\operatorname{id}_x:x\in\operatorname{Ob}(\mathcal{G})\}.
\]
There are only countably many such finite words. Consequently,
$\operatorname{Mor}(\mathcal{G})$ is at most countable.

It follows that every groupoid in $\mathbf{Gpd}_{f_0,f_g}$ is isomorphic to one whose object set is a finite ordinal
\[
[n]=\{1,\ldots,n\}
\]
and whose morphism set is a subset of $\mathbb{N}$. For fixed
$[n]$ and $M\subseteq\mathbb{N}$, a groupoid structure is determined
by the source, target, identity, inverse, and composition maps in the sets $[n]$ and $M$. There is only a set of possible choices for these maps.
Since the collections
\[
\{[n]:n\in\mathbb{N}\}
\qquad\text{and}\qquad
\mathcal{P}(\mathbb{N})
\]
are sets, there is a set of such coded groupoids. Therefore, there
exists a set of objects that contain an isomorphic representative of
every groupoid in $\mathbf{Gpd}_{f_0,f_g}$. The full subcategory on this set of representatives is small and equivalent to $\mathbf{Gpd}_{f_0,f_g}$. Thus, $\mathbf{Gpd}_{f_0,f_g}$ is essentially small.

We next define the monoidal product by disjoint union:
\[
\mathcal{G}\otimes\mathcal{H}
=
\mathcal{G}\sqcup\mathcal{H}.
\]
The object and morphism sets are
\[
\operatorname{Ob}(\mathcal{G}\sqcup\mathcal{H})
=
\operatorname{Ob}(\mathcal{G})
\sqcup
\operatorname{Ob}(\mathcal{H})
\ \
\mbox{ and}
 \ \
\operatorname{Mor}(\mathcal{G}\sqcup\mathcal{H})
=
\operatorname{Mor}(\mathcal{G})
\sqcup
\operatorname{Mor}(\mathcal{H}),
\]
with no morphisms between the two summands. If $S_{\mathcal{G}}$ and
$S_{\mathcal{H}}$ are finite generating sets for $\mathcal{G}$ and
$\mathcal{H}$, respectively, then
\[
S_{\mathcal{G}}\sqcup S_{\mathcal{H}}
\]
is a finite generating set for $\mathcal{G}\sqcup\mathcal{H}$. Hence, disjoint union preserves finite-object, finitely generated groupoids. For functors
\[
F\colon\mathcal{G}\to\mathcal{G}'
\qquad\text{and}\qquad
K\colon\mathcal{H}\to\mathcal{H}',
\]
their disjoint union defines a functor
\[
F\sqcup K\colon
\mathcal{G}\sqcup\mathcal{H}
\longrightarrow
\mathcal{G}'\sqcup\mathcal{H}'.
\]
Thus disjoint union defines a bifunctor
\[
\sqcup\colon\mathbf{Gpd}_{f_0,f_g}\times\mathbf{Gpd}_{f_0,f_g}\longrightarrow\mathbf{Gpd}_{f_0,f_g}.
\]
The empty groupoid $\varnothing$, which is generated by the empty set, is the monoidal unit. The canonical associativity, unit, and symmetry isomorphisms for disjoint unions provide the required natural isomorphisms, and their coherence follows from the coproduct construction \cite[Section~VII.1, p.~163; Section~XI.1]{maclane1998categories}. Therefore,
\[
(\mathbf{Gpd}_{f_0,f_g},\sqcup,\varnothing)
\]
is a symmetric monoidal category.

It remains to prove that $\mathbf{Gpd}_{f_0,f_g}$ has categorical images in the sense of \cite[Section~2]{patel2016generalized}. We first recall that a functor
\[
m\colon\mathcal{L}\to\mathcal{H}
\]
is a monomorphism in $\mathbf{Gpd}_{f_0,f_g}$ if and only if it is injective on objects and faithful. The reverse implication follows immediately from the definition of monomorphism. For the forward implication, injectivity on objects can be tested using the one object trivial groupoid, while faithfulness can be tested using the groupoid generated by one isomorphism between two objects. Both testing groupoids belong to $\mathbf{Gpd}_{f_0,f_g}$. Now, let
\[
F\colon\mathcal{G}\longrightarrow\mathcal{H}
\]
be a functor in $\mathbf{Gpd}_{f_0,f_g}$. Define $\operatorname{Im}(F)$ as the subgroupoid of $\mathcal{H}$ generated by the objects $F(x)$, for $x\in\operatorname{Ob}(\mathcal{G})$, and the morphisms $F(\alpha)$, for $\alpha\in\operatorname{Mor}(\mathcal{G})$. Thus,
\[
\operatorname{Ob}(\operatorname{Im}(F))
=
\{F(x):x\in\operatorname{Ob}(\mathcal{G})\}.
\]
If $S$ is a finite generating set for $\mathcal{G}$, then
$\operatorname{Im}(F)$ is generated by the finite set $F(S)$.
Therefore, $\operatorname{Im}(F)$ has finitely many objects, is
finitely generated, and therefore belongs to $\mathbf{Gpd}_{f_0,f_g}$. There is a factorization
\[
\mathcal{G}
\xrightarrow{\overline{F}}
\operatorname{Im}(F)
\xrightarrow{\iota_F}
\mathcal{H},
\]
where $\overline{F}$ is obtained by restricting the codomain of $F$ and $\iota_F$ is the inclusion functor. The functor $\iota_F$ is injective on objects and faithful, and is therefore a monomorphism.

To verify the universal property, suppose that $F$ admits another
factorization
\[
\mathcal{G}
\xrightarrow{Q}
\mathcal{L}
\xrightarrow{m}
\mathcal{H},
\qquad
F=m\circ Q,
\]
where $m$ is a monomorphism. Since $m$ is injective on objects, the
assignment
\[
F(x)\longmapsto Q(x)
\]
is well-defined. Since $m$ is faithful, the assignment
$F(\alpha)\longmapsto Q(\alpha)$
extends uniquely to a well-defined functor
\[
U\colon\operatorname{Im}(F)\longrightarrow\mathcal{L}.
\]
Indeed, if two generator composites represent the same morphism in
$\operatorname{Im}(F)$, then their images under $m$ are the same in
$\mathcal{H}$; faithfulness of $m$ implies that the same corresponding
composites in $\mathcal{L}$. By construction,
\[
m\circ U=\iota_F
\qquad\text{and}\qquad
U\circ\overline{F}=Q.
\]
The functor $U$ is unique because $\operatorname{Im}(F)$ is generated
by objects and morphisms in the image of $F$. Thus
, $\iota_F\colon\operatorname{Im}(F)\to\mathcal{H}$ is the categorical image of $F$.

Therefore, $\mathbf{Gpd}_{f_0,f_g}$ is an essentially small symmetric
monoidal category with categorical images, with a monoidal product given by the disjoint union and a unit given by the empty groupoid.
\end{proof}

The proof of our second main result is given below.

\begin{proof}[Proof of Theorem \ref{main2}]
Let $\Gamma$ be a finite simple undirected graph. By definition, the objects of its $A$-fundamental groupoid are the vertices of $\Gamma$:
\[
\operatorname{Ob}\bigl(\Pi_1^A(\Gamma)\bigr)=V(\Gamma).
\]
Since $\Gamma$ is finite, $\Pi_1^A(\Gamma)$ has finitely many objects. For vertices $x,y\in V(\Gamma)$, the morphisms from $x$ to $y$ are the path-homotopy classes of finite paths from $x$ to $y$:
\[
\operatorname{Hom}_{\Pi_1^A(\Gamma)}(x,y)
=
\bigl\{[\gamma]\mid
\gamma\colon x\rightsquigarrow y
\text{ is a finite path in }\Gamma
\bigr\}.
\]
The collection of all finite paths in $\Gamma$ is a set. Indeed, it is
contained in
\[
\coprod_{n\geq 0}V(\Gamma)^{n+1},
\]
which is countable because $V(\Gamma)$ is finite. Therefore, the collection of path-homotopy classes is also a set and $\Pi_1^A(\Gamma)$ is a small groupoid. It remains to prove that $\Pi_1^A(\Gamma)$ is finitely generated. Let \[\overrightarrow{E}(\Gamma)
=
\bigl\{(u,v)\in V(\Gamma)\times V(\Gamma)
\mid \{u,v\}\in E(\Gamma)\bigr\}\]
be the set of oriented edges of $\Gamma$. For every $(u,v)\in
\overrightarrow{E}(\Gamma)$, let
$ \varepsilon_{u,v}\colon I_1\longrightarrow \Gamma$
be the path of length one satisfying
$
\varepsilon_{u,v}(0)=u$ and $
\varepsilon_{u,v}(1)=v.$ Consider the set
\[
S_\Gamma=
\bigl\{[\varepsilon_{u,v}]
\mid (u,v)\in\overrightarrow{E}(\Gamma)\bigr\}
\subseteq
\operatorname{Mor}\bigl(\Pi_1^A(\Gamma)\bigr).
\]
Since $\Gamma$ has finitely many edges, $S_\Gamma$ is finite. We claim that $S_\Gamma$ generates $\Pi_1^A(\Gamma)$. Let
\[
\gamma\colon I_n\longrightarrow \Gamma
\]
be any finite path, and write
\[
x_i=\gamma(i),
\qquad 0\leq i\leq n.
\]
For every $1\leq i\leq n$, the definition of a graph map implies that
either
$
x_{i-1}=x_i
$
or
$
\{x_{i-1},x_i\}\in E(\Gamma).
$
Let $\gamma_i$ denote the path of length one from $x_{i-1}$ to $x_i$.
If $x_{i-1}=x_i$, then $[\gamma_i]$ is the identity morphism of
$x_i$. Otherwise,
\[
[\gamma_i]=[\varepsilon_{x_{i-1},x_i}]\in S_\Gamma.
\]
The original path is the concatenation of these elementary paths:
$\gamma=\gamma_1*\gamma_2*\cdots *\gamma_n.$
Since the composition in $\Pi_1^A(\Gamma)$ is induced by concatenation, it follows that
\[
[\gamma]
=
[\gamma_n]\circ\cdots\circ[\gamma_2]\circ[\gamma_1].
\]
Thus, every morphism in $\Pi_1^A(\Gamma)$ is a finite composite of elements of $S_\Gamma$ and identity morphisms. Equivalently, one may choose one orientation for every edge and use inverses to obtain the opposite orientations. Hence $\Pi_1^A(\Gamma)$ is generated finitely as a groupoid.

Finally, if
\[
f\colon \Gamma_1\longrightarrow \Gamma_2
\]
is a graph map between graphs, the induced functor
\[
\Pi_1^A(f)\colon\Pi_1^A(\Gamma_1)\longrightarrow\Pi_1^A(\Gamma_2)
\]
is defined on objects by
$x\longmapsto f(x)$
and on morphisms by $[\gamma]\longmapsto[f\circ\gamma].$
This is well-defined because graph maps preserve path homotopies and
 preserve identities, inverses, and concatenation. Since the
morphisms of $\mathbf{Gpd}_{f_0,f_g}$ are all functors, $\Pi_1^A(f)$
is a morphism in $\mathbf{Gpd}_{f_0,f_g}$. The identity and composition properties are derived from the functoriality of $\Pi_1^A$. Therefore,
\[
\Pi_1^A\colon
\mathbf{Graph}_{\mathrm{fin}}
\longrightarrow
\mathbf{Gpd}_{f_0,f_g}
\]
is a well-defined functor.
\end{proof}

Let $\mathcal{X}\colon(\mathbb{R},\leq)\longrightarrow \mathbf{Graph}_{\mathrm{fin}}$ be a filtration of finite simple undirected graphs. Since the $A$-homotopy groupoid defines a functor $\Pi_1^A\colon\mathbf{Graph}_{\mathrm{fin}} \longrightarrow\mathbf{Gpd}_{f_0,f_g}$, the composite
$$
F_{\mathcal{X}}
=
\Pi_1^A\circ\mathcal{X}\colon
(\mathbb{R},\leq)
\longrightarrow
\mathbf{Gpd}_{f_0,f_g}
$$
is a persistence module. Explicitly, $F_{\mathcal{X}}(p)=\Pi_1^A(\mathcal{X}(p))$ for every \(p\in\mathbb{R}\), and

$$
F_{\mathcal{X}}(p\leq q)
=
\Pi_1^A\bigl(\mathcal{X}(p\leq q)\bigr)
$$
for every \(p\leq q\). Suppose, moreover, that \(\mathcal{X}\) changes only at the values of a finite set $S=\{s_1<\cdots<s_n\}$,
 that $\mathcal{X}(p)=\varnothing$ for $p<s_1$, and that its structure maps are graph isomorphisms whenever $p$ and $q$ lie in the same interval determined by $S$. Since
$$
\Pi_1^A(\varnothing)=\varnothing
$$
and $\Pi_1^A$ sends graph isomorphisms to groupoid isomorphisms, $F_{\mathcal{X}}$ is $S$-constructible. Consequently,$$
F_{\mathcal{X}}
\in
\mathbf{PMod}\bigl(\mathbf{Gpd}_{f_0,f_g}\bigr).
$$
Since \(\mathbf{Gpd}_{f_0,f_g}\) is an essentially small symmetric monoidal category with categorical images, Patel's type construction $\mathcal A$ associates a generalized persistence diagram with $F_{\mathcal{X}}$.

\begin{example}\label{ex:C5-K5} \rm
Let \(C_5\) be the cycle graph on five vertices and let \(K_5\) be the
complete graph on the same vertex set. Consider the
\(\{1,2\}\)-constructible filtration
\[
\mathcal{X}\colon(\R,\leq)\longrightarrow\mathbf{Graph}_{\mathrm{fin}}
\]
defined by
\[
\mathcal{X}(t)=
\begin{cases}
\varnothing, & t<1,\\
C_5, & 1\leq t<2,\\
K_5, & 2\leq t,
\end{cases}
\]
whose nonidentity transition map is induced by the inclusion
\[
C_5\hookrightarrow K_5.
\]

\begin{figure}[ht]
\centering
\begin{tikzpicture}[
    scale = 0.8,
    vertex/.style={circle,fill=black,inner sep=2.2pt},
    edge/.style={line width=0.8pt},
    transition/.style={->,line width=0.9pt}
]

\node[font=\Large] (E) at (-5,0) {$\varnothing$};

\begin{scope}[shift={(0,0)}]
    \coordinate (C1) at (90:1.2);
    \coordinate (C2) at (18:1.2);
    \coordinate (C3) at (-54:1.2);
    \coordinate (C4) at (-126:1.2);
    \coordinate (C5) at (162:1.2);

    \draw[edge]
        (C1)--(C2)--(C3)--(C4)--(C5)--cycle;

    \foreach \i in {1,...,5}
        \node[vertex] at (C\i) {};
\end{scope}

\begin{scope}[shift={(5,0)}]
    \coordinate (K1) at (90:1.2);
    \coordinate (K2) at (18:1.2);
    \coordinate (K3) at (-54:1.2);
    \coordinate (K4) at (-126:1.2);
    \coordinate (K5) at (162:1.2);

    \draw[edge]
        (K1)--(K2)--(K3)--(K4)--(K5)--cycle;

    \draw[edge]
        (K1)--(K3)
        (K1)--(K4)
        (K2)--(K4)
        (K2)--(K5)
        (K3)--(K5);

    \foreach \i in {1,...,5}
        \node[vertex] at (K\i) {};
\end{scope}

\draw[transition]
    (-4.35,0)--(-1.55,0)
    node[midway,above] {$\iota_1$};

\draw[transition]
    (1.55,0)--(3.45,0)
    node[midway,above] {$\iota_2$};

\node at (0,1.65) {$C_5$};
\node at (5,1.65) {$K_5$};

\node at (-5,-1.65) {$t<1$};
\node at (0,-1.65) {$1\leq t<2$};
\node at (5,-1.65) {$2\leq t$};

\end{tikzpicture}
\caption{The $\{1,2\}$-constructible filtration
\(\varnothing\hookrightarrow C_5\hookrightarrow K_5\).}
\label{fig:C5-K5-filtration}
\end{figure}
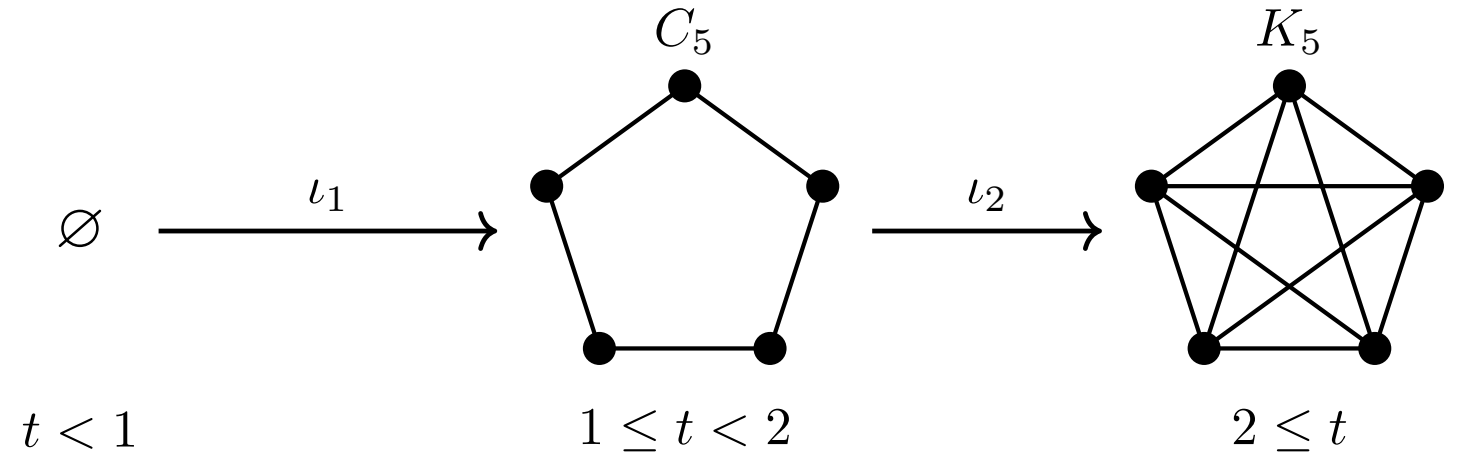
Applying the \(A\)-fundamental groupoid functor gives the persistence
module
\[
F=\Pi_1^A\circ\mathcal{X}\colon
(\R,\leq)\longrightarrow\mathbf{Gpd}_{f_0,f_g}.
\]

The \(A\)-fundamental groupoid of \(C_5\) has five objects. After choosing
an orientation of the cycle and reference paths between its vertices,
every morphism can be described by its endpoints together with an
integer recording its winding number. Consequently, the standard connected-groupoid decomposition \cite[Section~1.7, equation~(1.7.1)]{ronbrown2} gives a noncanonical isomorphism
\[
\Pi_1^A(C_5)
\cong
\operatorname{Pair}_5\times B\Z.
\]
On the other hand, every two \(A\)-paths in \(K_5\) having the same
endpoints are \(A\)-homotopic. Hence
\[
\Pi_1^A(K_5)
\cong
\operatorname{Pair}_5.
\]
Thus, up to the above identifications, the only nonidentity transition
in \(F\) is the functor
\[
\Phi\colon
\operatorname{Pair}_5\times B\Z
\longrightarrow
\operatorname{Pair}_5
\]
which is the identity on objects and forgets the winding number. More
precisely, it sends every morphism between a fixed pair of objects to
the unique morphism between those objects in
\(\operatorname{Pair}_5\). Therefore,
\[
\operatorname{im}(\Phi)=\operatorname{Pair}_5.
\]

Let
\[
\mathcal{G}
=
\operatorname{Pair}_5\times B\Z,
\qquad
\mathcal{H}
=
\operatorname{Pair}_5.
\]
The nonzero values of the image-based rank function that are relevant
for the Möbius inversion are
\[
dF_{K_0}([1,2))=[\mathcal{G}],
\qquad
dF_{K_0}([1,\infty))=[\mathcal{H}],
\qquad
dF_{K_0}([2,\infty))=[\mathcal{H}].
\]
Here the value on \([1,2)\) is computed, according to Patel's endpoint
convention, using a parameter immediately to the left of \(2\).

Choose \(s_0=0<1\). Since \(F(0)=\varnothing\), every structure map
starting at \(0\) has empty image, and hence
\[
dF_{K_0}([0,2))
=
dF_{K_0}([0,\infty))
=
[\varnothing]
=
0.
\]
Since \(2\) is the largest element of the constructibility set, the
M\"obius inversion formula gives
\[F_{K_0}([1,2))
=
dF_{K_0}([1,2))
-dF_{K_0}([1,\infty))
+dF_{K_0}([0,\infty))
-dF_{K_0}([0,2))
=
[\operatorname{Pair}_5\times B\Z]
-
[\operatorname{Pair}_5],
\]
and
\[F_{K_0}([1,\infty))
=
dF_{K_0}([1,\infty))
-dF_{K_0}([0,\infty))
=
[\operatorname{Pair}_5].\]
Finally,
\[F_{K_0}([2,\infty))
=
dF_{K_0}([2,\infty))
-dF_{K_0}([1,\infty))
=
[\operatorname{Pair}_5]-[\operatorname{Pair}_5]
=0.
\]
All other values of \(F_{K_0}\) are zero. Therefore, the generalized
persistence diagram is supported on the two intervals
\[
[1,2)
\qquad\text{and}\qquad
[1,\infty),
\]
with respective values
\[
[\operatorname{Pair}_5\times B\Z]
-
[\operatorname{Pair}_5]
\qquad\text{and}\qquad
[\operatorname{Pair}_5].
\]

The finite-interval term records the disappearance, at parameter \(2\),
of the nontrivial automorphism group generated by the cycle in \(C_5\).
The infinite-interval term records the groupoid structure that remains
after this cycle disappears. Notice that the difference
\[
[\operatorname{Pair}_5\times B\Z]
-
[\operatorname{Pair}_5]
\]
is a formal difference in the Grothendieck group
\(K_0(\mathbf{Gpd}_{f_0,f_g})\), rather than an actual difference of
groupoids.
\end{example}

\begin{example}[A filtration on five vertices] \rm
\label{ex:weighted-pentagon}
 Consider the graph filtration
\[
\mathcal{X}\colon(\R,\leq)\longrightarrow
\mathbf{Graph}_{\mathrm{fin}}
\]
specified in Figure~\ref{fig:weighted-pentagon-filtration}. The
filtration is constant between consecutive integer values and
\(\mathcal{X}(t)=\varnothing\) for \(t<0\).

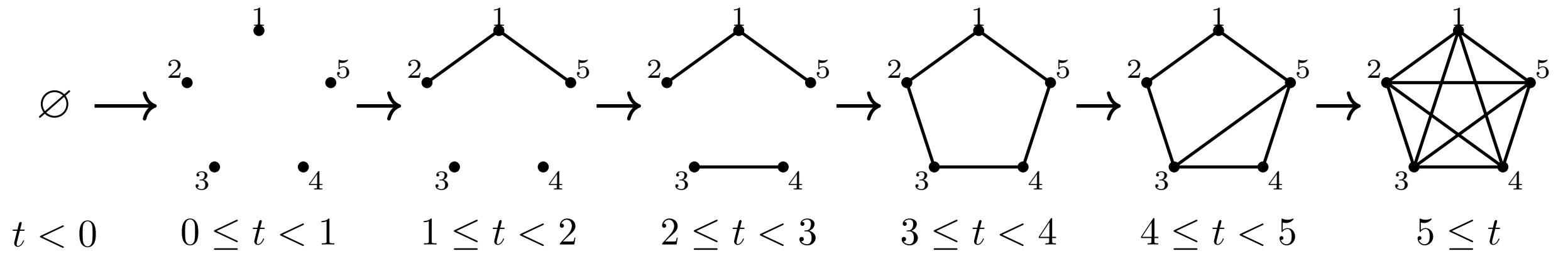
\begin{figure}[ht]
\centering
\resizebox{\textwidth}{!}{%
\begin{tikzpicture}[
    scale=0.8,
    vertex/.style={circle,fill=black,inner sep=1pt},
    vlabel/.style={font=\scriptsize,inner sep=1pt},
    edge/.style={line width=0.8pt},
    transition/.style={->,line width=0.9pt}
]

\node[font=\Large] (E) at (-2.3,0) {$\varnothing$};

\begin{scope}[shift={(0,0)}]
    \coordinate (A1) at (90:0.85);
    \coordinate (A2) at (162:0.85);
    \coordinate (A3) at (-126:0.85);
    \coordinate (A4) at (-54:0.85);
    \coordinate (A5) at (18:0.85);

    \foreach \i in {1,...,5}
        \node[vertex] at (A\i) {};

    \node[vlabel,above]       at (A1) {$1$};
    \node[vlabel,above left]  at (A2) {$2$};
    \node[vlabel,below left]  at (A3) {$3$};
    \node[vlabel,below right] at (A4) {$4$};
    \node[vlabel,above right] at (A5) {$5$};
\end{scope}

\begin{scope}[shift={(2.7,0)}]
    \coordinate (B1) at (90:0.85);
    \coordinate (B2) at (162:0.85);
    \coordinate (B3) at (-126:0.85);
    \coordinate (B4) at (-54:0.85);
    \coordinate (B5) at (18:0.85);

    \draw[edge]
        (B1)--(B2)
        (B1)--(B5);

    \foreach \i in {1,...,5}
        \node[vertex] at (B\i) {};

    \node[vlabel,above]       at (B1) {$1$};
    \node[vlabel,above left]  at (B2) {$2$};
    \node[vlabel,below left]  at (B3) {$3$};
    \node[vlabel,below right] at (B4) {$4$};
    \node[vlabel,above right] at (B5) {$5$};
\end{scope}

\begin{scope}[shift={(5.4,0)}]
    \coordinate (C1) at (90:0.85);
    \coordinate (C2) at (162:0.85);
    \coordinate (C3) at (-126:0.85);
    \coordinate (C4) at (-54:0.85);
    \coordinate (C5) at (18:0.85);

    \draw[edge]
        (C1)--(C2)
        (C1)--(C5)
        (C3)--(C4);

    \foreach \i in {1,...,5}
        \node[vertex] at (C\i) {};

    \node[vlabel,above]       at (C1) {$1$};
    \node[vlabel,above left]  at (C2) {$2$};
    \node[vlabel,below left]  at (C3) {$3$};
    \node[vlabel,below right] at (C4) {$4$};
    \node[vlabel,above right] at (C5) {$5$};
\end{scope}

\begin{scope}[shift={(8.1,0)}]
    \coordinate (D1) at (90:0.85);
    \coordinate (D2) at (162:0.85);
    \coordinate (D3) at (-126:0.85);
    \coordinate (D4) at (-54:0.85);
    \coordinate (D5) at (18:0.85);

    \draw[edge]
        (D1)--(D2)--(D3)--(D4)--(D5)--cycle;

    \foreach \i in {1,...,5}
        \node[vertex] at (D\i) {};

    \node[vlabel,above]       at (D1) {$1$};
    \node[vlabel,above left]  at (D2) {$2$};
    \node[vlabel,below left]  at (D3) {$3$};
    \node[vlabel,below right] at (D4) {$4$};
    \node[vlabel,above right] at (D5) {$5$};
\end{scope}

\begin{scope}[shift={(10.8,0)}]
    \coordinate (F1) at (90:0.85);
    \coordinate (F2) at (162:0.85);
    \coordinate (F3) at (-126:0.85);
    \coordinate (F4) at (-54:0.85);
    \coordinate (F5) at (18:0.85);

    \draw[edge]
        (F1)--(F2)--(F3)--(F4)--(F5)--cycle;

    \draw[edge]
        (F3)--(F5);

    \foreach \i in {1,...,5}
        \node[vertex] at (F\i) {};

    \node[vlabel,above]       at (F1) {$1$};
    \node[vlabel,above left]  at (F2) {$2$};
    \node[vlabel,below left]  at (F3) {$3$};
    \node[vlabel,below right] at (F4) {$4$};
    \node[vlabel,above right] at (F5) {$5$};
\end{scope}

\begin{scope}[shift={(13.5,0)}]
    \coordinate (K1) at (90:0.85);
    \coordinate (K2) at (162:0.85);
    \coordinate (K3) at (-126:0.85);
    \coordinate (K4) at (-54:0.85);
    \coordinate (K5) at (18:0.85);

    \draw[edge]
        (K1)--(K2)--(K3)--(K4)--(K5)--cycle;

    \draw[edge]
        (K1)--(K3)
        (K1)--(K4)
        (K2)--(K4)
        (K2)--(K5)
        (K3)--(K5);

    \foreach \i in {1,...,5}
        \node[vertex] at (K\i) {};

    \node[vlabel,above]       at (K1) {$1$};
    \node[vlabel,above left]  at (K2) {$2$};
    \node[vlabel,below left]  at (K3) {$3$};
    \node[vlabel,below right] at (K4) {$4$};
    \node[vlabel,above right] at (K5) {$5$};
\end{scope}

\draw[transition] (-1.85,0)--(-1.15,0);
\draw[transition] (1.10,0)--(1.60,0);
\draw[transition] (3.80,0)--(4.30,0);
\draw[transition] (6.50,0)--(7.00,0);
\draw[transition] (9.20,0)--(9.70,0);
\draw[transition] (11.90,0)--(12.40,0);

\node at (-2.3,-1.45) {$t<0$};
\node at (0,-1.45) {$0\leq t<1$};
\node at (2.7,-1.45) {$1\leq t<2$};
\node at (5.4,-1.45) {$2\leq t<3$};
\node at (8.1,-1.45) {$3\leq t<4$};
\node at (10.8,-1.45) {$4\leq t<5$};
\node at (13.5,-1.45) {$5\leq t$};

\end{tikzpicture}
}
\caption{A filtration on five vertices.}
\label{fig:weighted-pentagon-filtration}
\end{figure}

The filtration \(\mathcal{X}\) is constructible with respect to
\[
S_{\mathcal X}=\{0,1,2,3,4,5\}.
\]
Applying the \(A\)-fundamental groupoid functor gives a persistence
module
\[
F=\Pi_1^A\circ\mathcal{X}\colon
(\R,\leq)\longrightarrow\mathbf{Gpd}_{f_0,f_g}.
\]

Let \(\mathbf{1}\) denote the trivial one-object groupoid, and set
\[
\mathcal{D}
=
\bigsqcup_{i=1}^{5}\mathbf{1},\qquad \mathcal{T}_1
=
\operatorname{Pair}_3
\sqcup\mathbf{1}\sqcup\mathbf{1},\qquad
\mathcal{T}_2
=
\operatorname{Pair}_3
\sqcup\operatorname{Pair}_2,\qquad
\mathcal{C}
=
\operatorname{Pair}_5\times B\Z,\qquad
\mathcal{P}
=
\operatorname{Pair}_5.
\]
At time \(0\), the graph has five isolated vertices, so
\[
F(0)\cong\mathcal{D}.
\]
At time \(1\), the vertices \(1,2,5\) form a tree, while the vertices
\(3\) and \(4\) remain isolated. Hence
\[
F(1)\cong\mathcal{T}_1.
\]
At time \(2\), the graph has two tree components, with vertex sets
\(\{1,2,5\}\) and \(\{3,4\}\). Therefore,
\[
F(2)\cong\mathcal{T}_2.
\]
At time \(3\), the graph is \(C_5\), and hence
\[
F(3)
=
\Pi_1^A(C_5)
\cong
\mathcal{C}
=
\operatorname{Pair}_5\times B\Z.
\]

At time \(4\), the edge \(\{3,5\}\) produces the triangle
\[
3\longrightarrow4\longrightarrow5\longrightarrow3
\]
and the four-cycle
\[
1\longrightarrow2\longrightarrow3
\longrightarrow5\longrightarrow1.
\]
Both cycles are trivial under \(A\)-homotopy. Since they generate the
two ordinary cycles of the graph at this stage, every \(A\)-loop is
trivial. Consequently,
\[
F(4)\cong\mathcal{P}
=
\operatorname{Pair}_5.
\]
Finally,
\[
F(5)
=
\Pi_1^A(K_5)
\cong
\operatorname{Pair}_5
=
\mathcal{P}.
\]

The functor induced by the inclusion
\[
\mathcal{X}(4)\hookrightarrow\mathcal{X}(5)
\]
is an isomorphism at the level of \(A\)-fundamental groupoids.
Therefore, although the graph filtration changes at \(5\), the
persistence module \(F\) is already constructible with respect to
\[
S_F=\{0,1,2,3,4\}.
\]

Up to the preceding identifications, the persistence module has the
form
\[
\varnothing
\longrightarrow
\mathcal{D}
\longrightarrow
\mathcal{T}_1
\longrightarrow
\mathcal{T}_2
\longrightarrow
\mathcal{C}
\longrightarrow
\mathcal{P}
\xrightarrow{\cong}
\mathcal{P}.
\]

We now compute the image-based rank function. The inclusion-induced
functors are the identity on objects. Moreover, a morphism represented
by a path at an earlier stage remains represented by the same path at
all later stages. Thus, for
\(q\in\R\cup\{\infty\}\),
\begin{align*}
dF_{K_0}([0,q))
&=[\mathcal{D}],
&&q>0,\\
dF_{K_0}([1,q))
&=[\mathcal{T}_1],
&&q>1,\\
dF_{K_0}([2,q))
&=[\mathcal{T}_2],
&&q>2,\\
dF_{K_0}([3,q))
&=[\mathcal{C}],
&&3<q\leq4,\\
dF_{K_0}([3,q))
&=[\mathcal{P}],
&&q>4,\\
dF_{K_0}([4,q))
&=[\mathcal{P}],
&&q>4.
\end{align*}
Finite right endpoints in these expressions are interpreted using
Patel's left-limit convention.

The only nonzero finite-interval value of the Möbius inversion is
\[
F_{K_0}([3,4))
=
dF_{K_0}([3,4))
-dF_{K_0}([3,\infty))
+dF_{K_0}([2,\infty))
-dF_{K_0}([2,4))\\
=
[\mathcal{C}]
-[\mathcal{P}]
+[\mathcal{T}_2]
-[\mathcal{T}_2]\]
\[=
[\operatorname{Pair}_5\times B\Z]
-
[\operatorname{Pair}_5].
\]
The nonzero infinite-interval values are
\[F_{K_0}([0,\infty))
=
[\mathcal{D}]
=
5[\mathbf{1}],\ \ 
F_{K_0}([1,\infty))
=
[\mathcal{T}_1]-[\mathcal{D}]
=
[\operatorname{Pair}_3]-3[\mathbf{1}], \  \  
F_{K_0}([2,\infty))
=
[\mathcal{T}_2]-[\mathcal{T}_1]
=
[\operatorname{Pair}_2]-2[\mathbf{1}],\]
\[F_{K_0}([3,\infty))
=
[\mathcal{P}]-[\mathcal{T}_2]
=
[\operatorname{Pair}_5]
-[\operatorname{Pair}_3]
-[\operatorname{Pair}_2] \ \ \
F_{K_0}([4,\infty))
=
F_{K_0}([5,\infty))
=
0,
\]
and all other values of \(F_{K_0}\) are zero. The value
\[
[\operatorname{Pair}_5\times B\Z]
-
[\operatorname{Pair}_5]
\]
on the interval \([3,4)\) records the \(A\)-homotopy class created when
the five-cycle is completed at time \(3\) and killed when the edge
\(\{3,5\}\) is added at time \(4\). The infinite-interval values record
the vertices and the morphisms that connect them. Finally, adding the
remaining edges at time \(5\) produces no contribution to the
generalized persistence diagram, since it does not change the
\(A\)-fundamental groupoid up to isomorphism.
\end{example}

\section{The Argument of the Third Main Theorem}\label{s6}

We proceed to prove Theorem \ref{main3}, that is,  the  required structures---categorical images, essential smallness, and the symmetric monoidal condition---are satisfied under the given hypotheses.


\begin{proof}[Part 1 - Existence of Categorical Images - Proof of Theorem \ref{main3}]
In category theory, the image of a morphism $f: X \to Y$ is defined as the smallest subobject of $Y$ through which $f$ factors. Specifically, an image consists of an object $\operatorname{Im}(f)$, a monomorphism $m: \operatorname{Im}(f) \to Y$, and a morphism $e: X \to \operatorname{Im}(f)$ such that $f = m \circ e$. Furthermore, for any other factorization $f = m' \circ e'$ where $m': Z \to Y$ is a monomorphism, there must exist a unique morphism $h: \operatorname{Im}(f) \to Z$ such that $m' \circ h = m$ and $h \circ e = e'$.

Let $\mathcal{G}$ be a groupoid. By definition, every morphism $f \in \operatorname{Hom}_\mathcal{G}(X, Y)$ is an isomorphism. Because every isomorphism is simultaneously an epimorphism and a monomorphism, we can trivially factor $f$ as:
\begin{equation}
    f = \operatorname{id}_Y \circ f
\end{equation}
where $m = \operatorname{id}_Y$ and $e = f$. The identity morphism $\operatorname{id}_Y$ is clearly a monomorphism. Here $\operatorname{Im}(f) = Y$.

Now, suppose we have another factorization $f = m' \circ e'$ for some object $Z$ and monomorphism $m': Z \to Y$. Because $\mathcal{G}$ is a groupoid, the monomorphism $m'$ is necessarily an isomorphism. We construct the universal mapping diagram:

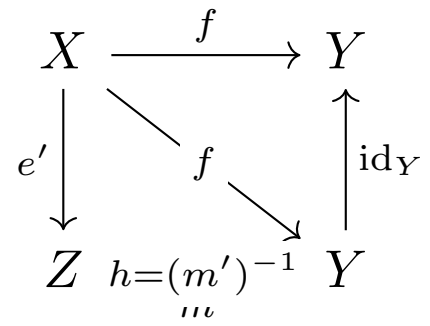
\begin{center}
\begin{tikzcd}[column sep=large, row sep=large]
X \arrow[r, "f"] \arrow[d, "e'"'] \arrow[dr, "f" description] & Y \\
Z \arrow[r, "m'"'] & Y \arrow[u, "\operatorname{id}_Y"'] \arrow[l, "h = (m')^{-1}" description, dashed]
\end{tikzcd}
\end{center}

To satisfy the universal property, we require a unique morphism $h: Y \to Z$ such that $m' \circ h = \operatorname{id}_Y$ and $h \circ f = e'$. Because $m'$ is an isomorphism, its inverse $(m')^{-1}$ exists and is unique. Setting $h = (m')^{-1}$ directly satisfies $m' \circ h = \operatorname{id}_Y$. Furthermore, 
\begin{equation}
    h \circ f = (m')^{-1} \circ (m' \circ e') = e'
\end{equation}
Thus, the condition is perfectly satisfied. The categorical image of any morphism $f: X \to Y$ in a groupoid is simply $(Y, \operatorname{id}_Y)$. Therefore, every groupoid automatically possesses categorical images.
\end{proof}


\begin{proof}[Part 2 - Essential Smallness - Proof of Theorem \ref{main3}]
A category is essentially small if it is equivalent to a small category; the following skeleton argument is standard \cite[Section~IV.4, pp.~93--94]{maclane1998categories}. 

Every category is equivalent to its skeleton. For a groupoid $\mathcal{G}$, the skeletal category $\mathcal{S}_{\mathcal{G}}$ is formed by choosing exactly one representative object for each connected component (isomorphism class) in $\mathcal{G}$. 
The objects of $\mathcal{S}_{\mathcal{G}}$ are indexed exactly by $\pi_0(\mathcal{G})$. 
Because the only morphisms in $\mathcal{S}_{\mathcal{G}}$ are the automorphisms of these representative objects, the collection of all morphisms in the skeleton is the disjoint union:
\begin{equation}
    \operatorname{Mor}(\mathcal{S}_{\mathcal{G}}) = \bigsqcup_{X \in \operatorname{Ob}(\mathcal{S}_{\mathcal{G}})} \operatorname{Aut}(X)
\end{equation}
If $\pi_0(\mathcal{G})$ is a set, and each vertex group $\operatorname{Aut}(X)$ is a set, then both $\operatorname{Ob}(\mathcal{S}_{\mathcal{G}})$ and $\operatorname{Mor}(\mathcal{S}_{\mathcal{G}})$ are sets. Therefore the skeleton $\mathcal S_{\mathcal G}$ is small. Conversely, an equivalence with a small category supplies a set of representatives and identifies each automorphism group with a set of morphisms in that category.
\end{proof}


\begin{proof}[Part 3 -  Eckmann-Hilton Argument - Proof of Theorem \ref{main3}]
Suppose we endow the groupoid $\mathcal{G}$ with a symmetric monoidal structure. This requires a bifunctor $\otimes: \mathcal{G} \times \mathcal{G} \to \mathcal{G}$ and a chosen unit object $I \in \operatorname{Ob}(\mathcal{G})$. 

Consider the vertex group of the monoidal unit, $\operatorname{Aut}(I) = \operatorname{Hom}_{\mathcal{G}}(I, I)$. This set inherits two distinct binary operations:
\begin{enumerate}
    \item Standard categorical composition: $\circ$
    \item The monoidal tensor product of morphisms: $\otimes$ (implicitly utilizing the unitor isomorphisms $\lambda_I = \rho_I: I \otimes I \xrightarrow{\sim} I$ to identify $\operatorname{Aut}(I \otimes I)$ with $\operatorname{Aut}(I)$).
\end{enumerate}
Because $\otimes$ is a bifunctor, it must satisfy the interchange law. For any four morphisms $f_1, f_2, g_1, g_2 \in \operatorname{Aut}(I)$:
\[    (f_2 \circ f_1) \otimes (g_2 \circ g_1) = (f_2 \otimes g_2) \circ (f_1 \otimes g_1)\]
Furthermore, the identity morphism $\operatorname{id}_I$ acts as the unit for both operations. Applying the classical Eckmann--Hilton Argument \cite[Theorem~1.3.1]{ronbrown2} (also \cite[Section~II.5, Exercise~5, p.~45]{maclane1998categories}), we evaluate the tensor product of two arbitrary automorphisms $f, g \in \operatorname{Aut}(I)$:
\[
    f \otimes g = (f \circ \operatorname{id}_I) \otimes (\operatorname{id}_I \circ g) \nonumber 
                = (f \otimes \operatorname{id}_I) \circ (\operatorname{id}_I \otimes g)  
                = f \circ g\]
Conversely, evaluating the elements in reverse order yields:
\[    f \otimes g = (\operatorname{id}_I \circ f) \otimes (g \circ \operatorname{id}_I) \nonumber 
                = (\operatorname{id}_I \otimes g) \circ (f \otimes \operatorname{id}_I) 
                = g \circ f\]
Equating the two results, we get:
\[    f \circ g = g \circ f\]
Therefore, $f \circ g = g \circ f$ for all $f, g \in \operatorname{Aut}(I)$. This proves that the endomorphism group of the monoidal unit in $\mathcal{G}$ must be strictly abelian. If $\mathcal{G}$ cannot house the monoidal unit in a connected component with an abelian fundamental group, the coherent symmetric tensor product cannot exist. Now, the symmetric monoidal structure was assumed throughout. The abelianness just proved is a necessary consequence, not a sufficient construction of that structure. The image and size assertions above complete the stated criterion.
\end{proof}

The notion of core of a category was probably introduced by MacLane and others for the first time.

\begin{definition}[See \cite{ronbrown1, ronbrown2, maclane1998categories}]
    Let $\mathbf{Grp}_{\mathrm{fg}}^{\simeq}$ be the \textit{core} of the category of finitely generated groups: its objects are finitely generated groups and its morphisms are group isomorphisms. 
\end{definition}

Taking all groups would not give an essentially small category. This observation is appropriate to show with details and we do it below.

\begin{proposition}\label{niceproposition}
$\mathbf{Grp}_{\mathrm{fg}}^{\simeq}$ is an essentially small symmetric monoidal category with images, with the tensor of the free product. 
\end{proposition}

For this proposition write $\mathcal C=\mathbf{Grp}_{\mathrm{fg}}^{\simeq}$. The coproduct construction of the tensor is classical \cite[Sections~VII.1 and XI.1]{maclane1998categories}. This core is an auxiliary example: its isomorphisms do not accommodate general transition maps in a graph filtration.

\begin{proof}[Proof of Proposition \ref{niceproposition}] We argue as done in the proof of Theorem \ref{main3}, showing analogies but also differences withe case of graphs.

\medskip

\textit{Part 1: Essential Smallness in $\mathbf{Grp}_{\mathrm{fg}}^{\simeq}$.}

\medskip

A category is essentially small if it is categorically equivalent to a small category (a category whose objects and morphisms both form mathematical sets, rather than proper classes). We must show that the isomorphism classes of $\mathcal{C}$ form a set, and the morphisms between any two objects form a set. Any finitely generated group with $n$ generators can be expressed as a quotient of the free group of rank $n$, denoted $F_n$. A quotient group is strictly defined by a normal subgroup of $F_n$. Because $F_n$ is a countably infinite group, the number of its subsets (and therefore the maximum possible number of normal subgroups) is bounded by the cardinality of the continuum, $2^{\aleph_0} = \mathfrak{c}$.
Taking the countable union over all possible finite ranks $n \in \mathbb{N}$, the total number of isomorphism classes of finitely generated groups is bounded by $\aleph_0 \times \mathfrak{c} = \mathfrak{c}$. Because $\mathfrak{c}$ is a well-defined mathematical set, the collection of isomorphism classes $\pi_0(\mathcal{C})$ forms a set. Furthermore, for any group $G \in \operatorname{Ob}(\mathcal{C})$, the morphisms from $G$ to itself are exactly its automorphism group, $\operatorname{Aut}(G)$. Because a finitely generated group is a countable set, the number of possible functions from $G$ to itself is bounded by $\mathfrak{c}$, which is also a set. Because the isomorphism classes form a set, and the morphism groups form sets, the skeleton of $\mathcal{C}$ is a small category. Therefore, $\mathcal{C}$ is essentially small.

\medskip

\textit{Part 2: Symmetric Monoidal Structure via Free Products in $\mathbf{Grp}_{\mathrm{fg}}^{\simeq}$.}

\medskip

We must show that the free product of groups endows $\mathcal{C}$ with a symmetric monoidal structure. This requires a bifunctor $\otimes$, a unit object $I$, and natural isomorphisms  satisfying Mac Lane's Coherence Conditions.

\medskip

\textit{Part 2.1. The Bifunctor and Unit:} We define the tensor product $G \otimes H$ as the free product of groups, $G \ast H$, and the monoidal unit $I$ as the trivial group $\{1\}$. To be a valid bifunctor in the core category, it must map isomorphisms to isomorphisms. Let $\phi: G_1 \xrightarrow{\sim} G_2$ and $\psi: H_1 \xrightarrow{\sim} H_2$ be group isomorphisms. By the universal property of the free product, there exists a unique group homomorphism $\phi \ast \psi: G_1 \ast H_1 \to G_2 \ast H_2$. Because both $\phi$ and $\psi$ are invertible, their free product extension is also invertible, meaning $\phi \ast \psi$ is an isomorphism and thus a valid morphism in $\mathcal{C}$.

\medskip

\textit{Part 2.2. Structural Isomorphisms:} In the parent category of groups ($\mathbf{Grp}$), the free product acts as the categorical coproduct. The universal property of the coproduct natively provides canonical homomorphisms for the associator ($\alpha$), left/right unitors ($\lambda, \rho$), and symmetry ($\sigma$). Crucially, all of these canonical maps are strictly group isomorphisms:   $ \alpha_{G,H,K}: (G \ast H) \ast K \xrightarrow{\sim} G \ast (H \ast K)$, \,
   $\lambda_G: \{1\} \ast G \xrightarrow{\sim} G$, \,
    $\rho_G: G \ast \{1\} \xrightarrow{\sim} G$ \,
    $\sigma_{G,H}: G \ast H \xrightarrow{\sim} H \ast G$. Because they are isomorphisms, they survive the restriction to the core and exist natively as valid morphisms in $\mathcal{C}$.

\medskip

\textit{Part 2.3. Coherence:} Any category with finite coproducts is automatically a symmetric monoidal category, meaning Mac Lane's pentagon and hexagon diagrams commute natively in $\mathbf{Grp}$. Because $\mathcal{C}$ is a subcategory of $\mathbf{Grp}$, any diagram of isomorphisms that commutes in $\mathbf{Grp}$ strictly commutes in $\mathcal{C}$. Thus, the free product endows $\mathcal{C}$ with a coherent symmetric monoidal structure.

We note explicitly that this construction successfully bypasses the Eckmann-Hilton obstruction. The endomorphism group of the monoidal unit in $\mathcal{C}$ is $\operatorname{Aut}(\{1\})$. Since the only automorphism of the trivial group is the identity map, this is the trivial group, which is abelian. Thus, the algebraic obstruction is resolved without restricting the non-abelian geometry of the other groups in the category.

\medskip

\textit{Part 3: Categorical Images in $\mathbf{Grp}_{\mathrm{fg}}^{\simeq}$.}

\medskip

In category theory, the image of a morphism $f: X \to Y$ is defined as the smallest subobject of $Y$ through which $f$ factors. Specifically, an image consists of an object $\operatorname{Im}(f)$ and a monomorphism $m: \operatorname{Im}(f) \to Y$ such that $f = m \circ e$ for some $e$. For any other factorization $f = m' \circ e'$ where $m': Z \to Y$ is a monomorphism, there must exist a unique morphism $h: \operatorname{Im}(f) \to Z$ such that $m' \circ h = m$ and $h \circ e = e'$.

Because $\mathcal{C}$ is a groupoid, every morphism $f \in \operatorname{Hom}_\mathcal{C}(X, Y)$ is an isomorphism. In category theory, every isomorphism is simultaneously an epimorphism and a monomorphism. We can trivially factor $f$ as:
\[    f = \operatorname{id}_Y \circ f \]
where the identity morphism $\operatorname{id}_Y$ is a monomorphism. Thus, our candidate image is $\operatorname{Im}(f) = Y$ with $m = \operatorname{id}_Y$ and $e = f$. Suppose that we have another factorization $f = m' \circ e'$ for some object $Z$ and a monomorphism $m': Z \to Y$. Because $\mathcal{C}$ is a groupoid, the monomorphism $m'$ is necessarily an isomorphism, meaning its inverse $(m')^{-1}$ exists. We construct the universal mapping diagram:

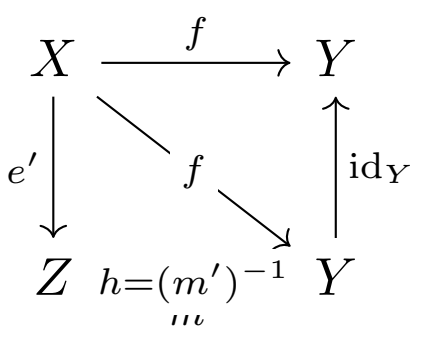
\begin{center}
\begin{tikzcd}[column sep=large, row sep=large]
X \arrow[r, "f"] \arrow[d, "e'"'] \arrow[dr, "f" description] & Y \\
Z \arrow[r, "m'"'] & Y \arrow[u, "\operatorname{id}_Y"'] \arrow[l, "h = (m')^{-1}" description, dashed]
\end{tikzcd}
\end{center}

We require a unique morphism $h: Y \to Z$ such that $m' \circ h = \operatorname{id}_Y$ and $h \circ f = e'$. Setting $h = (m')^{-1}$ directly satisfies $m' \circ h = \operatorname{id}_Y$. Furthermore, 
\[
    h \circ f = (m')^{-1} \circ (m' \circ e') = e'
\]
The universal condition is perfectly satisfied. The categorical image of any morphism $f: X \to Y$ in $\mathcal{C}$ is simply $(Y, \operatorname{id}_Y)$. Therefore, the category automatically possesses images.
\end{proof}

\section{A linear realization and type-\texorpdfstring{$\mathcal B$}{B} stability}
\label{sec:linearization}

Put $\mathcal C=\mathbf{Gpd}_{f_0,f_g}$ and fix a field $k$ of characteristic zero. Choose a small full subcategory $\mathcal D\subset\mathcal C$ containing one representative of each isomorphism class of nonempty groupoids. Such a choice is possible by Theorem~\ref{main1}. Throughout this section all functors between presheaves are natural transformations.

The Yoneda embedding and linear presheaves turn out to be relevant notions at this point: more precisely, the Yoneda embedding and its full faithfulness are classical notions which are illustrated in \cite[Section~III.2, pp.~61--62]{maclane1998categories}. For present purposes we use the following restriction and linearization:
\[
\mathcal P_k=[\mathcal D^{\mathrm{op}},\mathbf{Vect}_k],
\qquad
L_k(G)(H)=k[\operatorname{Hom}_{\mathcal C}(H,G)].
\]
Here $k[S]$ denotes the vector space freely generated by a set $S$, with $k[\varnothing]=0$. The restriction maps are the linear extensions of precomposition. A functor $f:G\to G'$ induces the natural transformation $L_k(f)$ by $[a]\mapsto[f\circ a]$. The category $\mathcal P_k$ is abelian, with kernels, cokernels and finite direct sums computed objectwise; this is the standard functor-category construction of \cite[Section~VIII.3, p.~199]{maclane1998categories}.

The elementary linearization of sets is classical (compare \cite[Sections~1.1--1.3]{deligne2004categorie} for finite sets). The closure below is a construction inside this fixed presheaf category. It is not the universal tensor-category construction in  \cite[Section 8]{deligne2004categorie}, and we claim no universal property for it.

\begin{lemma}\label{lem:normalized-faithful}
The functor $L_k:\mathcal C\to\mathcal P_k$ is faithful and sends the empty groupoid to the zero presheaf.
\end{lemma}
\begin{proof}
If $G\neq\varnothing$ and $f,g:G\to G'$ are distinct, choose $H\in\mathcal D$ and an isomorphism $a:H\to G$. Then $f\circ a\neq g\circ a$. The component at $H$ of $L_k(f)$ and of $L_k(g)$ sends the basis vector $[a]$ to two distinct basis vectors. Thus $L_k(f)\neq L_k(g)$. For an empty source there is at most one functor to a fixed target, so faithfulness is automatic. Finally, there is no functor from $H$ nonempty  to $\varnothing$,  so $L_k(\varnothing)(H)=0$.
\end{proof}

The restriction to nonempty test objects is necessary for the initial-value convention in persistence. If the empty test object were included, the unreduced presheaf would have value
\[
k[\operatorname{Hom}(\varnothing,\varnothing)]=k,
\]
so it would not send the empty groupoid to zero. Faithfulness does not mean fullness; arbitrary linear natural transformations need not arise from functors of groupoids. The abelian subcategory generated by the presheaves is useful in our discussions.

\begin{definition}\label{def:ak}
Let $\mathcal A_k$ be the smallest full subcategory of $\mathcal P_k$ containing all $L_k(G)$, closed under isomorphisms, finite direct sums  and kernels and cokernels of morphisms between its objects.
\end{definition}

We have all that we need to prove our final main result.

\begin{proof}[Proof of Theorem~\ref{main4}] We split the proof in two parts.

\medskip

\textit{Part 1 - Categorical Assertion.} 

\medskip

Closure under finite direct sums and fullness make $\mathcal A_k$ additive and $k$-linear. Its kernels and cokernels are those of the ambient abelian category. For each morphism the canonical map from its coimage to its image is therefore the ambient isomorphism. This proves that $\mathcal A_k$ is abelian; compare the axioms and image factorization in \cite[Section~VIII.3, pp.~198--200]{maclane1998categories}.

For essential smallness, start with a set of representatives of the objects $L_k(G)$, which exists by Theorem~\ref{main1}. At each stage adjoin finite direct sums and one kernel and one cokernel for every morphism between objects already chosen. There is only a set of such morphisms because $\mathcal D$ is small and $\mathcal P_k$ is locally small. The union of these stages over the nonnegative integers is a set. Its closure under isomorphisms is exactly $\mathcal A_k$: every pair of its objects occurs, up to isomorphism, at a common stage, and fullness ensures that the next stage supplies the required kernels and cokernels. Thus $\mathcal A_k$ is essentially small. There is also a direct concrete realization, namely
\[
U:\mathcal P_k\longrightarrow\mathbf{Vect}_k,
\qquad U(P)=\bigoplus_{H\in\operatorname{Ob}(\mathcal D)}P(H),
\qquad U(\alpha)=\bigoplus_H\alpha_H.
\]
This functor is faithful, since its restriction to each summand detects the corresponding component of a natural transformation. It is exact and preserves images: all these constructions are objectwise, and direct sums of vector spaces preserve exact sequences. Composing with the usual forgetful functor to sets gives the concrete-image realization required in \cite[Section~8, p.~413]{patel2016generalized}. Alternatively, the classical Freyd--Mitchell theorem supplies such a realization after passing to a small equivalent abelian category; see \cite[Introduction and Section~2]{montaruli}. We use the classical set-theoretic setting here, not the constructive strengthening investigated in that paper.
Finally, $\mathcal A_k$ is not equivalent to $\mathbf{Ab}$. Multiplication by $2$ on every object of the $k$-linear category $\mathcal A_k$ is invertible, since $\operatorname{char}k=0$. An equivalence of abelian categories preserves the addition on morphisms determined by biproducts. Such an equivalence would make multiplication by $2$ invertible on every abelian group, contrary to its action on $\mathbb Z$. The category $\mathcal A_k$ is nonzero: if $\mathbf1$ is the one-object trivial groupoid, $L_k(\mathbf1)(H)=k$ for every $H\in\mathcal D$. Lemma~\ref{lem:normalized-faithful} proves the remaining assertions about $L_k$.

\medskip

\textit{Part 2 - Constructibility and Continuous Stability}. 

\medskip

Let $M:(\mathbb R,\leq)\to\mathcal C$ be $S$-constructible. Then $\widehat M=L_k\circ M$ is constant up to isomorphism on the same parameter intervals, because functors preserve isomorphisms. Before the first critical value, $M$ is the empty groupoid, so Lemma~\ref{lem:normalized-faithful} gives $\widehat M=0$. Thus $\widehat M$ is constructible in the precise sense of \cite[Definition~2.2]{patel2016generalized}, for the direct-sum monoidal structure of $\mathcal A_k$.

Its type-$\mathcal B$ diagram is defined using \cite[Definition~6.2.1 and Definition~7.2]{patel2016generalized}. Given an actual $\epsilon$-interleaving of $M$ and $N$, postcomposition with $L_k$ preserves the two interleaving identities, so gives an $\epsilon$-interleaving of $\widehat M$ and $\widehat N$. Taking the infimum over all such shifts yields
\[
d_I(\widehat M,\widehat N)\leq d_I(M,N).
\]
No assumption that the infimum is attained is needed. Applying Theorem~\ref{stability}, which is \cite[Theorem~8.2]{patel2016generalized}, now gives
\[
d_E(\widehat M_{\mathcal B},\widehat N_{\mathcal B})
\leq d_I(\widehat M,\widehat N)
\leq d_I(M,N).
\]
\end{proof}

It is useful to note that the construction which we used in the proof of Theorem \ref{main4} changes the target category and then applies an existing stability theorem. It does not assert that $L_k$ preserves the categorical images in $\mathcal C$, or that the type-$\mathcal B$ diagram is the same invariant as the type-$\mathcal A$ diagram of the original groupoid module. Nor does faithfulness of $L_k$ establish nonvanishing or completeness of the resulting Grothendieck-group invariant. Those comparison and detection questions require separate study.

\bibliographystyle{plain}
\bibliography{references}

\end{document}